\documentclass[11pt]{article}
\usepackage[a4paper]{geometry}
\usepackage{amssymb, amsmath, amsthm, amsfonts, amscd}
\usepackage[english]{babel}
\usepackage{graphicx}
\usepackage{caption}
\usepackage{subcaption}
\usepackage{color,soul}
\usepackage{booktabs,dcolumn}                     
\usepackage[figureposition=bottom]{caption}       
\usepackage[usenames,dvipsnames,svgnames,table]{xcolor}
\usepackage{authblk}
\usepackage[title]{appendix}
\usepackage{tikz}
\usepackage{mathalfa}
\usepackage{comment}
\usepackage{accents}

\usetikzlibrary{arrows.meta, decorations.markings, angles, quotes}
\newtheorem{theorem}{Theorem}[section]

\newtheorem{lemma}[theorem]{Lemma}
\theoremstyle{definition}
\newtheorem{defn}[theorem]{Definition}

\newtheorem{prop}[theorem]{Proposition}
\newtheorem{remark}[theorem]{Remark}
\numberwithin{equation}{section}

\def\R{\mathbb{R}}
\def\C{\mathbb{C}}
\def\Re{\text{Re}}
\def\Im{\text{Im}}

\newcommand{\norm}[1]{\left\| #1 \right\|}
\newcommand{\abs}[1]{\left| #1 \right|}
\newcommand{\dd}{\mathrm{d}}

\title{Well--posedness of the {B}urgers' equation on the half--line through the {C}ole--{H}opf transformation}
\author{Huy Vo$^{\MakeLowercase{a}}$ and A. S. Fokas $^{\MakeLowercase{a,b}}$}
\affil{$^a$Department of Applied Mathematics and Theoretical Physics, University of Cambridge, UK}
\affil{$^b$Mathematics Research Center, Academy of Athens, Greece}

\begin{document}
\maketitle
\begin{abstract}
We rigorously establish the well--posedness of the initial-boundary value problem (IBVP) for the Burgers’ equation on the half-line for data in Sobolev-type spaces. The formal solution to this problem was previously derived via the Cole-Hopf transformation and the Fokas unified transform method under the assumption of smooth, Schwartz-class data. In this work, we prove the validity of this formal representation for initial data $u_0\in \tilde{H}^{s+1}_x (0,\infty)$ and boundary data $g_0 \in H_t^{{(2s+1)}/{4}}$ where $1/2<s<3/2$. well--posedness is established under two complementary conditions: (1) sufficiently small initial and boundary data, and (2) non-positive boundary data of arbitrary size. The latter condition utilizes the strong minimum principle for the corresponding linear heat equation to ensure the global invertibility of the nonlinear Cole-Hopf map. The proofs combine linear estimates for the heat IBVP with nonlinear bounds obtained through Taylor series representations and the algebraic properties of Sobolev spaces.
\end{abstract}

\section{Introduction}\label{sec:1}
The well--posedness of the Viscous Burgers' equation on the full line was establish by Daniel Dix in\cite{Dix}, utilising the dissipation effect of the heat equation to set up a contraction map. However, it is challenging to use this approach for the same problem on the half-line. The major obstacle is that for the heat equation on the half-line, due to the presence of the boundary data, the solution of the linear problem lacks the same dissipation effect. 

The Viscous Burgers' equation on the half-line is formally solved in \cite{Fokas_2014}. In this work, we explore the validity of the formula in section 5 of \cite{Fokas_2014} for data in Sobolev spaces. The above formula is derived assuming ``nice'' initial and boundary data (e.g.\ of Schwartz class). Is this formula valid for larger classes of data belonging to appropriate Sobolev-type spaces? In particular, we investigate the following initial-boundary value problem for the Burgers' equation formulated on the half-line:
\begin{subequations} \label{eq:1.1}
\begin{align}
& u_t = u_{xx}+ 2uu_x, && x \in (0, \infty),\ t \in (0, T),\ T < \infty, \label{eq:LS-eq} \\
&u(x, 0) = u_0(x) \in \tilde{H}^{s+1}_x(0, \infty), && x \in [0, \infty), \label{eq:LS-ic} \\
&u(0, t) = g_0(t) \in H^{(2s+1)/4}_t(0, T), && t \in [0, T]. \label{eq:LS-bc}
\end{align}
\end{subequations}

The Cole-Hopf transformation implies that 
\begin{equation}
u(x, t) = \frac{v_x(x, t)}{1 + v(x, t)}, \label{eq:1.2} 
\end{equation}
{where $v(x, t)$ solves the heat equation}
\begin{equation}
v_t(x, t) = v_{xx}(x, t) , \label{eq:1.3}
\end{equation}
{with the initial condition}
\begin{equation}
v_0(x) = e^{\int_0^x u_0(\xi) \, d\xi} - 1, \quad 0 < x < \infty, \label{eq:1.4}
\end{equation}
{and the Neumann boundary condition}
\begin{equation}
v_x(0, t) = h_1(t), \quad 0 < t < T. \quad \label{eq:1.5} 
\end{equation}
{The function $h_1(t)$ is determined as the unique solution of the
following linear Volterra integral equation:}
\begin{equation}
h_1(t) = F(t) - \frac{g_0(t)}{\sqrt{\pi}} \int_0^t \frac{h_1(s)}{\sqrt{t - s}} \, ds, \quad 0 < t < T, \label{eq:1.6}
\end{equation}
{with the forcing $F(t)$ being given in terms of $u_0(x)$ and
$g_0(t)$ by the formula}
\begin{equation}
F(t) = g_0(t) + \frac{g_0(t)}{\sqrt{\pi t}} \int_0^\infty e^{-\xi^2 / 4t} \left[ e^{\int_\infty^\xi u_0(\xi')\, d\xi'}  - 1 \right] d\xi. \label{eq:1.7}
\end{equation}
The solution $v(x, t)$ of the heat equation with the initial
condition $v_0(x)$ and the Neumann boundary condition $h_1(t)$ is given by
\begin{equation}
v(x, t) = \frac{1}{2\pi} \int_{-\infty}^\infty e^{i k x - k^2 t} \hat{v}_0(k) \, dk - \frac{1}{2\pi} \int_{\partial D^+} e^{i k x - k^2 t} \left[ \hat{v}_0(-k) + 2  \tilde{h}_1(k^2) \right] dk, \label{eq:1.8}
\end{equation}
\textit{where}
\begin{align}
\hat{v}_0(k) &= \int_0^\infty e^{-i k x} v_0(x) \, dx, \quad k \in \mathbb{C}, \quad \operatorname{Im} k \leqslant 0, \label{eq:1.9}  \\
\tilde{h}_1(k) &= \int_0^T e^{k s} h_1(s) \, ds, \quad k \in \mathbb{C}, \label{eq:1.10}
\end{align}
and $\partial D^+$ denotes the boundary of the region $D^+$ in the $k$-plane, depicted in Figure 1 and defined by $D^+=\{k\in\C:\Im(k)\geq |\Re (k)|\}$.

\begin{figure} \label{fig:1}
\center
\includegraphics[scale=2]{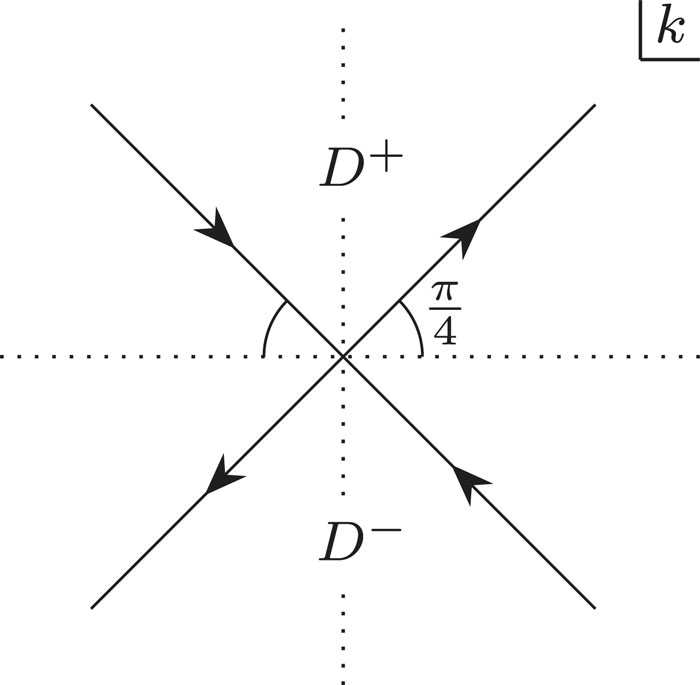}
\caption{The regions $D^{\pm}$ with their positively oriented boundaries $\partial D^{\pm}$.}
\end{figure}

\noindent Here, the spaces $H^s_x(0, \infty)$ and $H^{(2s+1)/4}_t(0, T)$ are defined as restrictions of the corresponding spaces over $\R$ according to the following definition.

\begin{defn}
For $s \in \R$ and $a, b \in \R$ (possibly infinite) with $b > a$, the Sobolev space $H^s(a, b)$ is defined by
\[
H^s(a, b) = \left\{ f : f = F|_{(a,b)} \text{ where } F \in H^s(\R) \right\},
\]
equipped with the norm
\[
\norm{f}_{H^s(a,b)} = \inf_{\substack{F \in H^s(\R) \\ F|_{(a,b)} = f}} \norm{F}_{H^s(\R)},
\]
where the Sobolev norm on $\R$ is defined via the Fourier transform of $F$ by
\[
\norm{F}_{H^s(\R)} = \left( \int_{\xi \in \R} (1 + \xi^2)^s \abs{\hat{F}(\xi)}^2 \, d\xi \right)^{1/2}.
\]
The space $\tilde{H}^{s+1}(0,\infty)$ is defined as 
\begin{align}
    \tilde{H}^{s+1}(0,\infty) &:= \left\{ f : \int^x_\infty f(\xi) d\xi  \in H^{s+1}(\R^+) \right\} \label{eq:1.11} \\
    &= \left\{ f \in H^{s}(\R^+) : \int^x_\infty f(\xi) d\xi  \in L^2(\R^+) \right\} . \label{eq:1.12}
\end{align}
\begin{equation}
    \norm{f}_{\tilde{H}^{s+1}(0,\infty)} := \norm{\int^x_\infty f(\xi) d\xi }_{H^{s+1}(\R^+)}. \label{eq:1.13}
\end{equation}
\end{defn}
One could think about the initial value $u_0$ as in $H^{s}(\R^+)$, supplemented with sufficient decay condition, a common condition for the Burgers' equation when analysed or evaluated using the Cole--Hopf transformation.

Clearly, there is an isomorphism between $ \tilde{H}^{s+1}(0,\infty)$ and ${H}^{s+1}(0,\infty)$. The set $\tilde{H}^{s+1}(0,\infty)$ can be estimated by smooth functions, but it is not an algebra.

For $s > 1/2$, the initial-boundary value problem ~\eqref{eq:LS-eq}--\eqref{eq:LS-bc} should be supplemented with the compatibility condition
\begin{equation}\label{eq:compatibility-1}
u_0(0) = g_0(0),
\end{equation}
which follows from the fact that $u_0(0)$ and $g_0(0)$ must both be equal to $u(0, 0)$.

Our results for the ibvp~\eqref{eq:LS-eq}--\eqref{eq:LS-bc} are as below.

\begin{theorem}[well--posedness of the ibvp~\eqref{eq:LS-eq}--\eqref{eq:LS-bc} \textbf{(Condition 1)}]\label{thm:1.1}
Suppose that $1/2<s< 3/2$.
Suppose further that 
\begin{align}
& \norm{g_0}_{\tilde{H}^{\frac{2s+1}{4}}_t (0,T)} \leq c_1, \label{eq:1.15}\\ 
   & \norm{u_0}_{\tilde{H}^{s+1}_x (0,\infty)} +\norm{g_0}_{\tilde{H}^{\frac{2s+1}{4}}_t(0,T)} \leq c_2, \label{eq:1.16}
\end{align}
where $c_1$ and $c_2$ sufficiently small are to be determined.
Then, the Cole-Hopf transformation formula \eqref{eq:1.2} defines a solution $u$ to the ibvp~\eqref{eq:LS-eq}--\eqref{eq:LS-bc} together with the compatibility condition~\eqref{eq:compatibility-1} in the following sense:

\begin{enumerate}
    \item \textbf{Space estimate:} The map $t \mapsto u(t)$ is $C([0, T]; \tilde{H}^{s+1}_x(0, \infty))$ with the estimate
    \[
    \sup_{t \in [0,T]} \norm{u(t)}_{\tilde{H}^{s+1}_x (0,\infty)} \le c_s \left( \norm{u_0}_{\tilde{H}^{s+1}_x (0,\infty)}  + \norm{g_0}_{H^{(2s+1)/4}_t(0,T)} \right).
    \]

    \item \textbf{Time estimates:} The map $x \mapsto u(x)$ is $C([0, \infty); H^{(2s+1)/4}_t(0, T))$ with the estimate
    \[
    \sup_{x \in [0,\infty)} \norm{u(x)}_{H^{(2s+1)/4}_t(0,T)} \le c_s \left( \norm{u_0}_{\tilde{H}^{s+1}_x (0,\infty)}  + \norm{g_0}_{H^{(2s+1)/4}_t(0,T)} \right).
    \]

    \item \textbf{Weak solution:} For any $\varphi(x) \in C_0^\infty(\R)$ such that $\varphi(0) = 0$, the function $\langle u(t), \varphi \rangle_{L^2_x(0,\infty)}$ is differentiable in $(0, T)$ with
    \[
    \frac{d}{dt} \langle u(t), \varphi \rangle_{L^2_x(0,\infty)} = \langle u(t), \varphi_{xx} \rangle_{L^2_x(0,\infty)} + g_0(t) \varphi_x(0) - \langle u^2(t),\varphi_x\rangle.
    \]

    \item \textbf{Initial and boundary conditions:} $u(x, 0) = u_0(x)$ for all $x \in [0, \infty)$, and $u(0, t) = g_0(t)$ for all $t \in [0, T]$.
\end{enumerate}

In the above estimates, $c_s > 0$ is a constant that depends on $s$.
\end{theorem}

\begin{theorem}[well--posedness of the ibvp~\eqref{eq:LS-eq}--\eqref{eq:LS-bc} \textbf{(Condition 2)}]\label{thm:1.2}
Suppose that $1/2<s< 3/2$.
Suppose further that 
\begin{align}
& \norm{g_0}_{\tilde{H}^{\frac{2s+1}{4}}_t (0,T)} \leq c_3, \label{eq:1.17}\\ 
   & \norm{u_0}_{\tilde{H}^{s+1}_x (0,\infty)} +\norm{g_0}_{\tilde{H}^{\frac{2s+1}{4}}_t(0,T)} \leq c_4, \label{eq:1.18}
\end{align}
where $c_3$ and $c_4$ can be arbitrarily large. Further, suppose that $g_0 \leq 0$ for all $t\in [0,T]$.
Then, the Cole-Hopf transformation formula \eqref{eq:1.2} defines a solution $u$ to the ibvp~\eqref{eq:LS-eq}--\eqref{eq:LS-bc} together with the compatibility condition~\eqref{eq:compatibility-1} in the following sense:

\begin{enumerate}
    \item  The map $t \mapsto u(t)$ is $C([0, T]; \tilde{H}^{s+1}_x(0, \infty))$.

    \item  The map $x \mapsto u(x)$ is $C([0, \infty); H^{(2s+1)/4}_t(0, T))$.

    \item  For any $\varphi(x) \in C_0^\infty(\R)$ such that $\varphi(0) = 0$, the function $\langle u(t), \varphi \rangle_{L^2_x(0,\infty)}$ is differentiable in $(0, T)$ with
    \[
    \frac{d}{dt} \langle u(t), \varphi \rangle_{L^2_x(0,\infty)} = \langle u(t), \varphi_{xx} \rangle_{L^2_x(0,\infty)} + g_0(t) \varphi_x(0) - \langle u^2(t),\varphi_x\rangle.
    \]

    \item $u(x, 0) = u_0(x)$ for all $x \in [0, \infty)$, and $u(0, t) = g_0(t)$ for all $t \in [0, T]$.
\end{enumerate}
The exact upper bound for $\sup_{t \in [0,T]} \norm{u(t)}_{\tilde{H}^{s+1}_x (0,\infty)} $ and $\sup_{x \in [0,\infty)} \norm{u(x)}_{H^{(2s+1)/4}_t(0,T)}$ can be given explicitly in term of $\norm{u_0}_{\tilde{H}^{s+1}_x (\R^+)}$ and $\norm{g_0}_{{H}^{\frac{2s+1}{4}}_t ((0,T))}$, but the exact formula is rather complicated. However, it can be retrieved by following along sections \ref{sec:6}, \ref{sec:7}, and \ref{sec:8}.

In the above estimates, $c_s > 0$ is a constant that depends on $s$.
\end{theorem}

The strategy of this paper is as follow: given initial data $u_0$ and $g_0$, we first show that the boundary data $v_0$ and $h_1$ of the heat equation, using the Cole-Hopf transformation, are in Sobolev spaces, and depend continuously on $u_0$ and $g_0$. From standard linear estimates, we have that $v(x,t)$ satisfies certain space and time estimates. Then the inverse Cole-Hopf implies space and time estimate for $u(x,t)$. Although the transformations are nonlinear, it is possible to estimate the relevant maps using their Taylor series form, together with the fact that $H^s$ is an algebra whenever $s>1/2$ (Kato--Ponce's lemma). 

In analysing the inverse Cole-Hopf transform
\begin{equation*}
    u(x,t)= \frac{v_x(x,t)}{1+v(x,t)},
\end{equation*}
to ensure that this map is well defined, we need $1+v(x,t)>0$ in the domain $(0,\infty)\times (0,T]$. This is not an issue if the IBVP \eqref{eq:1.1} is posed on the full real line, since by the weak minimum principle of the heat equation, $1+v(x,t)\geq \min\{1+v_0(x)\}>0$. On the half-line, conditions 1 and 2 implies two different approaches to fulfill this requirement: for condition 1, $v(x,t)$ will be sufficiently small, so that $1+v(x,t)>0$; for condition 2, by using a version of the strong minimum principle of the heat equation (Lemma 6.3), we have $1+v(x,t)\geq \min\{1+v_0(x)\}>0$.

We define the constant $c_a$ as follow.
\begin{defn}
    Let $c_a>0$ be a sufficiently large constant such that the following algebraic properties hold for the Sobolev spaces
    \begin{itemize}
        \item For any $f$, $g\in H^s(\R^+)$, $\norm{fg}_{H^s(\R^+)} \leq c_a \norm{f}_{H^s(\R^+)} \norm{g}_{H^s(\R^+)} .$
    \item For any $f$, $g\in H^{s+1}(\R^+)$, $\norm{fg}_{H^{s+1}(\R^+)} \leq c_a \norm{f}_{H^{s+1}(\R^+)} \norm{g}_{H^{s+1}(\R^+)} .$
    \item For any $f$, $g\in H^{\frac{2s+1}{4}}_t(0,2)$, $\norm{fg}_{ H^{\frac{2s+1}{4}}_t(0,2)} \leq c_a \norm{f}_{ H^{\frac{2s+1}{4}}_t(0,2)} \norm{g}_{ H^{\frac{2s+1}{4}}_t(0,2)} .$
    \item For any $f$, $g\in H^{\frac{2s-1}{4}}_t(0,2)$, $\norm{fg}_{ H^{\frac{2s-1}{4}}_t(0,2)} \leq c_a \norm{f}_{ H^{\frac{2s-1}{4}}_t(0,2)} \norm{g}_{ H^{\frac{2s-1}{4}}_t(0,2)} .$
    \item For any $f$, $g\in H^{\frac{2s+1}{4}}_t(0,T)$, $\norm{fg}_{ H^{\frac{2s+1}{4}}_t(0,T)} \leq c_a \norm{f}_{ H^{\frac{2s+1}{4}}_t(0,T)} \norm{g}_{ H^{\frac{2s+1}{4}}_t(0,T)} .$
    \item For any $f$, $g\in H^{\frac{2s-1}{4}}_t(0,T)$, $\norm{fg}_{ H^{\frac{2s-1}{4}}_t(0,T)} \leq c_a \norm{f}_{ H^{\frac{2s-1}{4}}_t(0,T)} \norm{g}_{ H^{\frac{2s-1}{4}}_t(0,T)} .$
    \end{itemize}
\end{defn}

\section{Linear Estimates for the Heat IBVP}\label{sec:2}
We summarize here the needed rigorous results for the heat equation with Neumann boundary condition on the half-line. Let $s'=s+1$ to denote the desired Sobolev regularity of the heat equation. Then $\frac{3}{2}<s'<\frac{5}{2}.$

Consider the following ibvp for the heat equation:
\begin{subequations}
\label{eq:forced-linear-IBVP}
\begin{align}
&v_t - v_{xx}= 0  && x \in (0,\infty),\ t \in (0, T),\ T < 1, \label{eq:forced-linear-IBVP-a} \\
&v(x,0) = v_0(x) \in H^{s'}_x(0,\infty), && x \in [0,\infty), \label{eq:forced-linear-IBVP-b} \\
&v_x(0, t) = h_1(t) \in H^{\frac{2{s'}-1}{4}}_t(0,T), && t \in [0, T]. \label{eq:forced-linear-IBVP-c}
\end{align}
\end{subequations}
We require the following compatibility condition between the boundary data
\begin{equation}
    v_{0,x}(0)=h_1(0). \label{eq:2.2}
\end{equation}

The solution of the ibvp ~\eqref{eq:forced-linear-IBVP} is given explicitly by using the Fokas method:
\begin{multline}
v(x, t) = S[v_0, h_1](x, t)
= \frac{1}{2\pi} \int_{k \in \mathbb{R}} e^{ikx - k^2 t} \widehat{v}_0(k)\, dk
+ \frac{1}{2\pi} \int_{k \in \partial D^+} e^{ikx - k^2 t} \widehat{v}_0(-k)\, dk \\
\quad - \frac{1}{\pi} \int_{k \in \partial D^+} e^{ikx - k^2 t} \widetilde{h}_1(k^2, T)\, dk, \label{eq:UTM-formula}
\end{multline}
where
\begin{align}
\widehat{v}_0(k) &= \int_0^\infty e^{-ikx} v_0(x)\, dx, \quad\label{eq:half-line-FT} \\
\widetilde{h}_1(k^2, T) &= \int_0^T e^{ik^2 t} h_1(t)\, dt, \quad k \in \mathbb{C}, \label{eq:time-transform-g1}
\end{align}
and the contour of integration $\partial D^+$ is the positively oriented boundary of the first quadrant of the complex $k$-plane, as shown in Fig 1.

Using \eqref{eq:UTM-formula}, it is possible to establish the well--posedness result given below.

\begin{theorem}[well--posedness of the heat IBVP]\label{thm:2.1}
\label{thm:linear-wellposedness}
Suppose that $\frac{3}{2} < s' < \frac{5}{2}$. Then equation \eqref{eq:UTM-formula} defines a solution $v = S[v_0, h_1]$ of the IBVP~\eqref{eq:forced-linear-IBVP} supplemented with the compatibility condition (2.2), which admits the following space and time estimates:
\begin{align}
\sup_{t \in [0,T]} \|v(t)\|_{H^{s'}_x(0,\infty)} 
&\leq c_{s'} \Bigl( \|v_0\|_{H^{s'}_x(0,\infty)} + \|h_1\|_{H^{\frac{2s'-1}{4}}_t(0,T)}  \Bigr), \label{eq:linear-space-est} \\
\sup_{x \in [0,\infty)} \|v(x)\|_{H^{\frac{2s'+1}{4}}_t(0,T)} 
&\leq c_{s'} \Bigl( \|v_0\|_{H^{s'}_x(0,\infty)} + \|h_1\|_{H^{\frac{2s'-1}{4}}_t(0,T)}  \Bigr), \\
\sup_{x \in [0,\infty)} \|v_x(x)\|_{H^{\frac{2s'-1}{4}}_t(0,T)} 
&\leq c_{s'} \Bigl( \|v_0\|_{H^{s'}_x(0,\infty)} + \|h_1\|_{H^{\frac{2s'-1}{4}}_t(0,T)}  \Bigr), \label{eq:linear-time-est}
\end{align}

where $c_{s'} = c(s') > 0$ is a constant depending only on $s'$.
\end{theorem}

For the purpose of establishing the linear estimates of Theorem \ref{thm:2.1}, it is convenient to express the linear IBVP as a superposition of simpler problems.

\subsection*{Linear decompositions}
We decompose the IBVP \eqref{eq:forced-linear-IBVP} as follows.

\medskip

\noindent Let $V_0 \in H^{s'}(\mathbb{R})$ be the extension to the real line of $v_0$, satisfying
\begin{equation}
    \|V_0\|_{H^{s'}_x(\mathbb{R})} \leq c \|v_0\|_{H^{s'}_x(0,\infty)},
    \label{eq:extension-estimate}
\end{equation}
we express the IBVP as the superposition of the following two problems:

\begin{enumerate}
    \item The homogeneous linear IVP:
    \begin{subequations}
        \label{eq:IVP-A1}
        \begin{align}
            V_t &= V_{xx} , && x \in \mathbb{R},\ t \in (0,T), \label{eq:IVP-A1-a} \\
            V(x,0) &= V_0(x) \in H^{s'}_x(\mathbb{R}), && x \in \mathbb{R}, \label{eq:IVP-A1-b}
        \end{align}
    \end{subequations}
    with solution
    \begin{equation}
        V(x,t) = S[V_0](x,t) = \frac{1}{2\pi} \int_{\xi \in \mathbb{R}} e^{i\xi x - \xi^2 t} \widehat{V}_0(\xi)\, d\xi,
        \label{eq:sol-IVP-A1}
    \end{equation}
    where $\widehat{V}_0$ is the Fourier transform of $V_0$ on the line given by
    \begin{equation}
        \widehat{V}_0(\xi) = \int_{x \in \mathbb{R}} e^{-i\xi x} V_0(x)\, dx, \quad \xi \in \mathbb{R}.
        \label{eq:FT-U0}
    \end{equation}
    \item The homogeneous linear IBVP with zero initial datum:
    \begin{subequations}
        \label{eq:IBVP-A2}
        \begin{align}
           v_t  &= v_{xx}, && x \in (0,\infty),\ t \in (0,T), \label{eq:IBVP-A2-a} \\
            v(x,0) &= 0, && x \in [0,\infty), \label{eq:IBVP-A2-b} \\
            v_x(0,t) &= {h^\dagger}(t) := h_1(t) - V_x(0,t), && t \in [0,T], \label{eq:IBVP-A2-c}
        \end{align}
    \end{subequations}
    with solution $v = S[0, {h^\dagger}]$ given by equation ~\eqref{eq:UTM-formula} after setting $v_0 = 0$.
\end{enumerate}
\begin{theorem}[Homogeneous linear IVP estimates]\label{thm:2.2}
\label{thm:homog-IVP-est}
The solution $V = S[V_0]$ of the linear IVP~\eqref{eq:IVP-A1} given by formula~\eqref{eq:sol-IVP-A1} admits the following estimates:
\begin{align}
\sup_{t \in [0,T]} \|V(t)\|_{H^{s'}_x(\mathbb{R})} &= \|V_0\|_{H^{s'}_x(\mathbb{R})}, && {s'} \in \mathbb{R}, \label{eq:IVP-space-est} \\
\sup_{x \in \mathbb{R}} \|V(x)\|_{H^{\frac{2s'+1}{4}}_t(0,T)} &\leq c_{s'} \|V_0\|_{H^{s'}_x(\mathbb{R})}, && s' \geq -\frac{1}{2}, \label{eq:IVP-time-est-U} \\
\sup_{x \in \mathbb{R}} \|V_x(x)\|_{H^{\frac{2s'-1}{4}}_t(0,T)} &\leq c_{s'} \|V_0\|_{H^{s'}_x(\mathbb{R})}, && s' \geq \frac{1}{2}. \label{eq:IVP-time-est-Ux}
\end{align}
\end{theorem}
\begin{theorem}[Linear estimates for compactly supported boundary data]\label{thm:2.3}
\label{thm:IBVP-compact-data}
The solution $v = S[0, {h^\dagger}]$ of the reduced IBVP \eqref{eq:IBVP-A2} given by equation \eqref{eq:UTM-formula} with $v_0=0$ admits the following space and time estimates
\begin{align}
\sup_{t \in [0,2]} \|v(t)\|_{H^{s'}_x(0,\infty)} &\leq c_{s'} \|{h^\dagger}\|_{H^{\frac{2{s'}-1}{4}}_t(\mathbb{R})}, && {s'} \geq \frac{1}{2}, \label{eq:IBVP-space-est} \\
\sup_{x \in [0,\infty)} \|v_x(x)\|_{H^{\frac{2{s'}-1}{4}}_t(0,2)} &\leq c_{s'} \|{h^\dagger}\|_{H^{\frac{2{s'}-1}{4}}_t(\mathbb{R})}, && {s'} \in \mathbb{R}, \label{eq:2.18} \\
\sup_{x \in [0,\infty)} \|v(x)\|_{H^{\frac{2{s'}+1}{4}}_t(0,2)} &\leq c_{s'} \|{h^\dagger}\|_{H^{\frac{2{s'}-1}{4}}_t(\mathbb{R})}, && s \geq \frac{1}{2}.
\label{eq:IBVP-time-est}
\end{align}
\end{theorem}
The proofs of the theorems of this section are presented in appendix A. The proofs follow the strategy of the works \cite{heat,NLS,NLS_Neumann}. However, special care must be taken when there is singularity in the frequency space. Also, the estimate \eqref{eq:IBVP-time-est} has not appeared in earlier work, and the proof here works only for ${s'}>\frac{1}{2}$.

Further, we noted the below result.

\begin{lemma} \label{lem:2.4}
    The solution $v=S[v_0,h_1]$ in Theorem \ref{thm:2.1} is a classical real solution. That is, assume $v_0$ and $h_1$ are real--valued functions, $v(x,t)$ is real for all $(x,t)\in[0,\infty)\times[0,T]$, and $v_t=v_{xx}$ for $(x,t)\in(0,\infty)\times(0,T]$ in the classical sense.
\end{lemma}
\begin{proof}
    See Appendix A.
\end{proof}

\section{The Cole-Hopf Transformation: Initial Data}\label{sec:3}
Formally, the initial condition $v_0$ is $e^{\int^x_\infty u_0(\xi) d\xi} -1$. We now define $v_0$
using the Taylor representation. The function $v_0$ is the convergence in $H^{s+1}(0,\infty)$ of the series 
\begin{equation}
    v_0:=\int^x_\infty u_0(\xi) d\xi +\frac{\left (\int^x_\infty u_0(\xi) d\xi \right)^2}{2!}+\frac{\left (\int^x_\infty u_0(\xi) d\xi \right)^3}{3!}+... \label{eq:3.1} 
\end{equation}
This series converges in $H^{s+1}(0,\infty)$, from the fact that $H^{s+1}(0,\infty)$ is an algebra. The radius of convergence for the above series is $\infty$, provided that $\norm{\int^x_\infty u_0(\xi) d\xi}_{H^{s+1}(0,\infty)} \leq c_2$. Further, given the uniform bound $c_2$, $\norm{v_0}_{H^{s+1}(0,\infty)}$ depends continously on $\norm{u_0}_{\tilde{H}^{s+1}}$.

Indeed, given $u_{0,1}$, $u_{0,2} \in \tilde{H}^{s+1}_{\R^+}$, then
\begin{multline}
 {v_{0,1}-v_{0,2} }  = \bigg( \int^x_\infty u_{0,1}(\xi) d\xi - \int^x_\infty u_{0,2}(\xi) d\xi \bigg) 
    \bigg( 1+ \frac{\int^x_\infty u_{0,1}(\xi) d\xi + \int^x_\infty u_{0,2}(\xi) d\xi }{ 2!} +\\ \ \frac{\left( \int^x_\infty u_{0,1}(\xi) d\xi \right)^2+ \int^x_\infty u_{0,1}(\xi) d\xi \int^x_\infty u_{0,2}(\xi) d\xi  + \left( \int^x_\infty u_{0,2}(\xi) d\xi\right) }{ 3!} +... \bigg). \label{eq:3.2} 
\end{multline}
Thus, using the uniform bound $c_2$ on $\norm{u_{0,1}}_{\tilde{H}^{s+1}}$ and $\norm{u_{0,2}}_{\tilde{H}^{s+1}}$, it follows that
\begin{multline}
    \norm{v_{0,1}-v_{0,2}}_{H^{s'}(\R^+)} \leq \norm{u_{0,1}-u_{0,2}}_{ \tilde{H}^{s+1}_{\R^+}} \left(1+c_2+\frac{c_2^2c_a}{2!}+\frac{c^3_2c^2_a}{3!}+... \right)\\ \leq \norm{u_{0,1}-u_{0,2}}_{ \tilde{H}^{s+1}_{\R^+}} \left( \frac{e^{c_2c_a} -1}{c_a}+1 \right). \label{eq:3.3}
\end{multline}
This shows that $v_0$ depends continuously on $u_0$, in the respective norms.
Estimating $\int^x_\infty u_0(\xi) d\xi$ using smooth functions, it follows that at $x=0$, the trace of $v_0$ is $e^{\int^0_\infty u(\xi) d\xi} -1$.

\section{The Cole-Hopf Transformation: Boundary Data}\label{sec:4}
For $T<1$, and for the initial condition $v_0 \in {H^{s'}(0,\infty)}$ and the Neumann boundary condition $h_1 \in H^{\frac{2s'-1}{4}}_t (0,T)$, the decomposition strategy used in linear analysis, expresses the solution in the form
\begin{equation}
    v(x,t)= S[V_0]+S[0,h_1-V_x]. \label{eq:4.1}
\end{equation}
Extending $(h_1-V_x)$ to $(0,2)$ and using a cut-off function $\theta$ as explained below, we can view $S[0,h_1-V_x]$ as the restriction of $S[0,(h_1-V_x)\theta(t)]$ to $(0,T)$. Thus, the unified transform formula gives 
\begin{equation}
    v(x,t)= \frac{1}{2\pi}\int_\R e^{ikx-k^2t} \hat{V_0}(k) dk -\frac{1}{\pi}\int_{\partial \mathcal{D}^+} e^{ikx-k^2t} \int_0^2 e^{k^2 s} (h_1(s)- V_x(0,s) ) \theta(s)\ ds . \label{eq:4.2}
\end{equation}
where:
\begin{enumerate}
    \item $\theta(s)$ is a smooth cut-off function, such that $\theta(s)$ is supported on $(-2,2)$, and $\theta \equiv 1$ on $(0,T)$,
    \item $v(x,t)$ is the restriction to $(0,\infty)\times (0,T)$ of the above formula.
\end{enumerate}
Then both components of the above sum satisfy the space and time estimates listed in section 2. Furthermore, they depend continuously on the initial and boundary data.
Following the approach used in section 5 of \cite{Fokas_2014}, we will employ the equation
\begin{equation}
    v_x(0,t)=g_0(t)[1+v(0,t)]. \label{eq:4.3}
\end{equation}
Formally 
\begin{multline}
    (h_1(t)-V_x(0,t))\theta(t)=\theta(t) \bigg[g_0(t)+\frac{g_0(t)}{2\pi} \int_\R e^{-k^2t} \hat{V_0}(k) dk - V_x(0,t) \bigg] \\-\frac{g_0(t)\theta(t)}{\pi} \int_{\partial \mathcal{D}^+} e^{-k^2t} \int_0^2 e^{k^2 s}( h_1(s)- V_x(0,s) ) \theta(s)\ ds\ dk. \label{eq:4.4}
\end{multline}
Consider the second integral of \eqref{eq:4.4}. For sufficiently smooth $h_1(t)$ and $V_0$ under the compatibility condition \eqref{eq:2.2}, using integration by part, 
\[
\int_0^2 e^{k^2 s}( h_1(s)- V_x(0,s) ) \theta(s) ds
\]
has $\mathcal{O}(\frac{1}{k^4})$ decay rates as $|k|\rightarrow\infty$. Thus, the outer integral converges absolutely. Furthermore,
\[
\int_{\partial \mathcal{D}^+} e^{-k^2(t-s)}dk 
\]
converges weakly along $\partial \mathcal{D}^+$. A limiting argument, using the additional kernel $e^{-\epsilon|k|^2}$, justify the rearrangement of the second integral of \eqref{eq:4.4}:
\begin{multline}
    (h_1(t)-V_x(0,t))\theta(t)=\theta(t) \bigg[g_0(t)+\frac{g_0(t)}{2\pi} \int_\R e^{-k^2t} \hat{V_0}(k) dk - V_x(0,t) \bigg] \\-\frac{g_0(t)\theta(t)}{\pi} \int_0^2  (h_1(s)- V_x(0,s) ) \theta(s)\int_{\partial \mathcal{D}^+} e^{-k^2(t-s)}dk  ds \\
    =\theta(t) \bigg[g_0(t)+\frac{g_0(t)}{2\pi} \int_\R e^{-k^2t} \hat{V_0}(k) dk - V_x(0,t) \bigg] \\-\frac{g_0(t)\theta(t)}{\pi} \int_0^t  (h_1(s)- V_x(0,s) ) \theta(s)\frac{1}{\sqrt{t-s}}  ds. \label{eq:4.5} 
\end{multline}
Let 
\begin{equation}
    F(t):=\theta(t) \bigg[g_0(t)+\frac{g_0(t)}{2\pi} \int_\R e^{-k^2t} \hat{V_0}(k) dk - V_x(0,t) \bigg], \label{eq:4.6}
\end{equation}
\begin{equation}
    {h^\ddagger}(t)=(h_1(t)-V_x(0,t))\theta(t), \label{eq:4.7}
\end{equation}
\eqref{eq:4.5} can be written in the form
\begin{equation}
    {h^\ddagger}(t)
    = F(t)-\frac{g_0(t)\theta(t)}{\pi} \int_0^t  {h^\ddagger}(s)\frac{1}{\sqrt{t-s}}  ds. \label{eq:4.8} \\
\end{equation}

Solving the above Volterra equation, we can retrieve the function $h_1$ as the restriction of ${h^\ddagger}(t) +V_x(t) $ to $(0,T)$.

\begin{lemma}
   The function ${h^\ddagger}(t)$ on $(0,2)$ as the solution of the Volterra equation (4.8) in $H^{\frac{2s'-1}{4}}_{0,t} (0,2)$ satisfies the followings:
\begin{enumerate}
    \item ${h^\ddagger}(t)$ exists and is unique.
    \item ${h^\ddagger}(t)$ depends continuously on $\norm{u_0}_{H^{s'}_x}$ and $\norm{g_0}_{H^{\frac{2s+1}{4}}_t}$.
\end{enumerate}
\end{lemma}

\begin{lemma} \label{lem:4.2}
        The function $h_1(t)$ defined above satisfies the followings:
\begin{enumerate}
    \item $h_1(t)$ is uniquely determined in $H^{\frac{2s'-1}{4}}_t (0,T)$.
    Further, $h_1(0)=v_x(0)$, thus, it satisfies the compatibility condition for the Neumann heat equation.
    \item When $u_0$ and $g_0$ are smooth, $h_1(t)= g_0(t)[1+v(0,t)]$, with $v(0,t)$ given by the unified transform formula.
    \item $h_1(t)$ depends continuously on $\norm{u_0}_{H^{s'}_x}$ and $\norm{g_0}_{H^{\frac{2s+1}{4}}_t}$.
\end{enumerate}
\end{lemma}
We first prove that there exists a unique solution to the Volterra equation.

\begin{prop}\label{prop:4.3}
    Given $f(t) \in H^{\frac{2s'-1}{4}}_{0,t} (0,2)$, the map
   \[
   R: f(t)\mapsto-\frac{g_0(t)\theta(t)}{\pi} \int_0^t  f(s)\frac{1}{\sqrt{t-s}} ds ,
   \] maps $f$ to the function $Rf \in H^{\frac{2s'-1}{4}}_{0,t} (0,2)$. 
    Furthermore,
\[
\norm{Rf}_{H^{\frac{2s'-1}{4}}_{t} (0,2)} \leq c \norm{g_0}_{\frac{2s+1}{4}} \norm{f}_{H^{\frac{2s'-1}{4}}_{t} (0,2)}.
\]  
\end{prop}
\begin{proof}
    Extending $f$ by zero outside of $(0,2)$, it follows that $\norm{f}_{H^{\frac{2s'-1}{4}}_{t}(\R)} = \norm{f}_{H^{\frac{2s'-1}{4}}_{t}(0,2)} $. Define the function
    \begin{equation}
        I(t) = 
\bigg\{
    \begin{array}{lr}
        \frac{1}{\sqrt{t}}, &\text{if } t\in (0,2), \\
        0, &\text{if } t\in \R\setminus(0,2).
    \end{array} \label{eq:4.9}
    \end{equation}
 $I(t)\in L^1(\R)$, with $\norm{I(t)}_{L^1(\R)}=2\sqrt{2}$. $\int_0^t  f(s)\frac{1}{\sqrt{t-s}} ds $ is the restriction to $(0,2)$ of the following $ {H^{\frac{2s'-1}{4}}_{t}(\R)}$ function $(f*I)(t)$. Thus,
\begin{multline}
    \norm{\int_0^t  f(s)\frac{1}{\sqrt{t-s}} ds}_{H^{\frac{2s'-1}{4}}_{t}(0,2) } \leq \norm{f*I}_{H^{\frac{2s'-1}{4}}_{t}(\R) } \leq \norm{f}_{H^{\frac{2s'-1}{4}}_{t}(\R)} \norm{I(t)}_{L^1(\R)} \\ \leq c \norm{f}_{H^{\frac{2s'-1}{4}}_{t}(0,2)}. \label{eq:4.10}
\end{multline}
For $\frac{1}{2}<\frac{2s'-1}{4}< 1$, since $\frac{1}{\sqrt{t-s}}$ is weakly integrable at zero and $\norm{f}_{L^\infty} < \infty$, 
\begin{equation}
 \lim_{t\rightarrow0}   \int_{0}^{t}  f(s)\frac{1}{\sqrt{t-s}} ds = 0. \label{eq:4.11}
\end{equation}
The trace of $\int_t^0  f(s)\frac{1}{\sqrt{t-s}} ds$ at $t=0$ is $0$. We conclude that $Rf(0)=0$.
Also, given that $g$ is a function in $H^{\frac{2s+1}{4}}_{t}(0,2) $, $\theta(t)$ vanishes at $2$, $Rf(2)=0$. Thus, $Rf \in H^{\frac{2s'-1}{4}}_{0,t} (0,2)$. The estimate of $Rf$ comes from the algebra property of $H^{\frac{2s+1}{4}}_{t}(0,2) $.
\end{proof}
Hence, for sufficiently small $g_0(t)$, the map $R$ is a contraction in $ H^{\frac{2s'-1}{4}}_{0,t} (0,2)$.
\begin{prop}\label{prop:4.4}
    $F(t)$ is a function in ${H^{\frac{2s'-1}{4}}_{0,t}(0,2) }$.
\end{prop}
\begin{proof}

    From standard linear estimate results, $F$ is a function in ${H^{\frac{2s'-1}{4}}_{t}(0,2) }$. For $\frac{1}{2}<\frac{2s'-1}{4}<\frac{3}{2}$, we need to check the traces at $0$ and $2$. 
    Clearly, by the cut-off function $\theta$, $F(2)=0$. To compute the trace at $0$, when $u_0$ is a smooth function, we use the following:
    \begin{align*}
        F(0)&=g_0(0)+g_0(0)v(0)-v_{0,x}(0) \\
        &= g_0(0)+g_0(0)(e^{\int^0_\infty u(\xi) d\xi} -1)-\partial_x (e^{\int^0_\infty u(\xi) d\xi} -1) \\
        &= g_0(0)+g_0(0)(e^{\int^0_\infty u(\xi) d\xi} -1) - u(0) e^{\int^0_\infty u(\xi) d\xi} \\
        &= 0 \quad \text{(using the compatibility condition $u_0(0)=g_0(0)$ at the origin.)}
    \end{align*}
\end{proof}
We are now ready to prove that the Volterra equation to find ${h^\ddagger}(t)$ is well-posed in $ H^{\frac{2s'-1}{4}}_{0,t} (0,2)$. Since $T$ is a contracting linear map, the solution, if exists, must be unique.
Choose ${h^\ddagger}(t)$ to be the convergence point of the following series in $ H^{\frac{2s'-1}{4}}_{0,t} (0,2)$:
\begin{equation}
    {h^\ddagger}(t)= (I+R+R^2+...)[F(t)]. \label{eq:4.12}
\end{equation}
Clearly ${h^\ddagger}(t)$ satisfies the Volterra equation. In this way, we derive at the unique solution. Furthermore, from standard linear estimate, we have: 
\begin{multline}
    \norm{F(t)}_{ H^{\frac{2s'-1}{4}}_{0,t} (0,2)} =\norm{\theta(t) \bigg[g_0(t)+\frac{g_0(t)}{2\pi} \int_\R e^{-k^2t} \hat{V_0}(k) dk - V_x(0,t) \bigg] }_{ H^{\frac{2s'-1}{4}}_{0,t} (0,2)} \\ \leq c_s\bigg( \norm{g_0}_{H^{\frac{2s+1}{4}}_{t} (0,2) } +\norm{g_0}_{H^{\frac{2s+1}{4}}_{t} (0,2) }\norm{v_0}_{H^s(\R^+)}+\norm{v_0}_{H^s(\R^+)} \bigg ). \label{eq:4.13}
\end{multline}
Thus,
\begin{multline}
    \norm{{h^\ddagger}(t) }_{ H^{\frac{2s'-1}{4}}_{0,t} (0,2)} \leq \frac{1}{1-\norm{R}_{\text{op}}}c_s\bigg( \norm{g_0}_{H^{\frac{2s+1}{4}}_{t} (0,2) } +\norm{g_0}_{H^{\frac{2s+1}{4}}_{t} (0,2) }\norm{v_0}_{H^s(\R^+)} +\norm{v_0}_{H^s(\R^+)} \bigg )\\
    \leq \frac{1}{1-c\norm{g_0}_{H^{\frac{2s+1}{4}}_{t} (0,2)}} c_s\bigg( \norm{g_0}_{H^{\frac{2s+1}{4}}_{t} (0,2) } +\norm{g_0}_{H^{\frac{2s+1}{4}}_{t} (0,2) }\norm{v_0}_{H^s(\R^+)} +\norm{v_0}_{H^s(\R^+)} \bigg ). \label{eq:4.14}
\end{multline}
To prove that ${h^\ddagger}(t)$ depends continuously on $g_0$ and $u_0$, we will use again the expression \eqref{eq:4.12}.
\begin{prop}\label{prop:4.5}
    The function ${h^\ddagger}(t)$, $t\in (0,2)$, depends continously on $u_0$ and $g_0$.
\end{prop}
\begin{proof}
    Let $(u_{0,m},g_{0,m})$, $(u_{0,n},g_{0,n})$ be two pairs of initial and boundary conditions of the Burgers' equation. From Section \ref{sec:2}, $\norm{v_{0,m}-v_{0,n}}_{H^{s'}_x(\R^+)} \lesssim \norm{u_{0,m}-u_{0,n}}_{\tilde{H}^{s+1}_x(\R^+)}$. Denote by $(F_m, R_m)$ the pair corresponding to $(u_{0,m},g_{0,m})$, and by $(F_n, R_n)$ the pair corresponding to $(u_{0,n},g_{0,n})$. Standard linear IBVP analysis implies
    \begin{equation}
        \norm{F_m(t)-F_n(t)}_{H^{\frac{2s'-1}{4}}_{0,t}(0,2)} \lesssim \norm{u_{0,m}-u_{0,n} }_{\tilde{H}^{s+1}_x(\R^+) }+\norm{g_{0,m}-g_{0,n} }_{H^{\frac{2s+1}{4}}_{t}(\R^+)}. \label{eq:4.15}
    \end{equation}
    Also, 
    \[
    (R_m-R_n)(F)(t)= -\frac{1}{\pi}(g_{0,m}-g_{0,n})\theta(t)  \int_0^t  F(s)\frac{1}{\sqrt{t-s}}\  ds.
    \]
    Thus,
    \begin{equation}
        \norm{R_m-R_n}_{op} \lesssim \norm{g_{0,m}-g_{0,n} }_{H^{\frac{2s+1}{4}}_{t}(0,2)} \lesssim \norm{g_{0,m}-g_{0,n} }_{H^{\frac{2s+1}{4}}_{t}(\R^+)}. \label{eq:4.16}
    \end{equation}
    We introduce the following notation: 
    \begin{align*}
        \norm{F_m(t)-F_n(t)}_{H^{\frac{2s'-1}{4}}_{t}(0,2)} &= \epsilon,\\
         \norm{R_m-R_n}_{op} &= \delta,\\
         \sup \left \{\norm{F_m(t)}_{H^{\frac{2s'-1}{4}}_{0,t}},\norm{F_n(t)}_{H^{\frac{2s'-1}{4}}_{0,t}} \right \}&=a,\\
          \sup \{\norm{R_m}_{op},\norm{R_n}_{op}\}&=b<1,
    \end{align*}
    where the last condition comes from contraction mapping assumption. We have
    \begin{align}
        {h^\ddagger}_m(t) &= (I+R_m+R^2_m+...)(F_m), \label{eq:4.17} \\
        {h^\ddagger}_n(t) &= (I+R_n+R^2_n+...)(F_n). \label{eq:4.18} 
    \end{align}
    Since both series converges absolutely in ${H^{\frac{2s'-1}{4}}_{0,t}(0,2)}$, we can arrange $ {h^\ddagger}_m(t)- {h^\ddagger}_n(t) $ in the form
    \begin{multline}
        (F_m-F_n)+(R_mF_m-R_nF_n)+(R_m^2F_m-R_n^2F_n)+(R^3_mF_m-R^3_n F_n)+... \label{eq:4.19}
    \end{multline}
    Now, we have
    \begin{align*}
        \norm{F_m-F_n} &\leq \epsilon.\\
        \norm{R_mF_m-R_nF_n} &\leq \norm{(R_m-R_n)F_m}+\norm{R_n(F_m-F_n)} \\
        &\leq \delta a +  b\epsilon.\\
        \norm{R_m^2F_m-R_n^2F_n} &\leq \norm{(R_m-R_n)R_nF_n}+\norm{R_n(R_nF_n-R_mF_m)} \\
        &\leq \delta b a + b(\delta a+ b\epsilon) = 2\delta ab+ b^2 \epsilon. \\
        ...\\
        \norm{R_m^{(k)} F_m -R_n^{(k)}F_n} &\leq k \delta a b^{k-1} + b^k \epsilon.\\
        ...
    \end{align*}
    Since $b<1$, the infinite sums below converge, and we have
    \begin{multline}
        \norm{{h^\ddagger}_m(t)-{h^\ddagger}_n(t) }_{ H^{\frac{2s'-1}{4}}_{0,t} (0,2)} \leq \delta a (1+2b+3b^2+4b^3+...)+ \epsilon(1+b+b^2+b^3+...) \\
        \lesssim \delta+\epsilon \lesssim  \norm{u_{0,m}-u_{0,n} }_{\tilde{H}^{s+1}_x(\R^+) }+\norm{g_{0,m}-g_{0,n} }_{H^{\frac{2s+1}{4}}_{t}(\R^+)}. \label{eq:4.20}
    \end{multline}
    Thus, we have proven that ${h^\ddagger}_0(t)$ depends continuously on $u_0$ and $g_0$.
\end{proof}
\begin{proof}[Proof of Lemma 4.2]
    $h_1(t)= {h^\ddagger}(t)+V_x(t)$. $V_x(0,t)$ is well-understood from the estimate on linear analysis problems, and since the trace at zero of ${h^\ddagger}(t)$ is zero, it follows that
\begin{equation}
    h_1(0)=V_x(0,0)= v_{0,x}(0). \label{eq:4.21}
\end{equation}
Thus, we have the compatibility condition for the heat equation.
Continuous dependence of $h_1(t)$ on $u_0$ and $g_0$ comes from Lemma 1 and the linear estimate of $V_x(t).$ 
\end{proof}

\begin{remark}
    We could still derive the existence and uniqueness of ${h^\ddagger}(t)$ without the bound
    \[
    \norm{g_0}_{H^{\frac{2s+1}{4}}_t(\R^+)} <c_2,
    \]
    for instance, by multiplying by $e^{-Mt}$ both sides of the Volterra equation. However, we still require the boundary data to be small in later sections, which requires $g_0$ to be sufficiently small. Thus, for simplicity, we assumed that $g_0$ is already sufficiently small.
\end{remark}
\section{well--posedness under Condition 1}\label{sec:5}
Recall that $u(x,t) = \frac{v_x(x,t)}{1+v(x,t)}=\partial_x \ln (1+v)$. For the rigorous analysis, we will use the Taylor series representation, 
\begin{equation}
    u(x,t)= \partial_x ( v-\frac{v^2}{2}+\frac{v^3}{3}-\frac{v^4}{4}+...), \label{eq:5.1}
\end{equation}
for space estimates, and 
\begin{equation}
    u(x,t)= ( v_x-vv_x+v^2v_x-...)=v_x(1-v+v^2-...), \label{eq:5.2}
\end{equation}
for time estimate. We need the second representation to avoid taking partial derivative in the $x$ direction, when we fix $x$ and let $t$ vary. Note that since $(0,T)$ is bounded, $1\in H^{\frac{2s'-1}{4}}_{t} (0,T)$. These representations agree pointwise when the boundary data are smooth.
\subsection{Space estimate}\label{subsec:5.1}

Combining all the above estimates, $v$ is bounded by $u_0$ and $g_0$ uniformly both in space and time. Hence, for sufficiently small upper bounds $c_1$ of $\norm{u_0}_{\tilde{H}^s_x}$ and $c_2$ of $\norm{g_0}_{H^{\frac{2s+1}{4}}_t}$, $\norm{v(\cdot,t)}_{H^{s}_x(\R^+)} \leq c<\frac{1}{c_a}$. Hence, the series 
\[
v-\frac{v^2}{2}+\frac{v^3}{3}-\frac{v^4}{4}+...
\]
converges to an element $q(x,t)$ in $H^{s'}$. $u(x,t)=\partial_x q(x,t)$ is in $H_x^{s'-1}$ or $H^s_x$, and is continuously dependent on $u_0$ and $g_0$. Further,
\begin{equation}
    \int^x_\infty u(\xi,t)d\xi = q(\xi,t)\bigg|^x_\infty = q(x,t), \label{eq:5.3}
\end{equation}
the antiderivative exists classically, since $s>1/2$. Thus, $u(\cdot,t)\in \tilde{H}^{s+1}_x(\R^+)$. To prove that $u(\cdot,t)$ depends continuously on $u_0$ and $g_0$, it suffices to prove that $u(\cdot,t)$ depends continuously on $v(\cdot,t)$. Note that:
\begin{align}
    \norm{u(\cdot,t)}_{\tilde{H}^{s+1}_x(\R^+) } &= \norm{\int^x_\infty u(\xi,t)d\xi}_{{H}^{s+1}_x(\R^+) } \label{eq:5.4}\\
    &\leq \norm{v(\cdot,t)}_{{H}^{s+1}_x(\R^+) } + \frac{c_a \norm{v(\cdot,t)}_{{H}^{s+1}_x(\R^+) }^2 }{2} + \frac{c_a^2 \norm{v(\cdot,t)}_{{H}^{s+1}_x(\R^+) }^3 }{3}+... \label{eq:5.5}\\
    &\leq \norm{v(\cdot,t)}_{{H}^{s+1}_x(\R^+) } + {c_a \norm{v(\cdot,t)}_{{H}^{s+1}_x(\R^+) }^2 } + {c_a^2 \norm{v(\cdot,t)}_{{H}^{s+1}_x(\R^+) }^3 }+... \label{eq:5.6}\\
    &\leq \frac{\norm{v(\cdot,t)}_{{H}^{s+1}_x(\R^+) } }{1-c_a \norm{v(\cdot,t)}_{{H}^{s+1}_x(\R^+) }}. \label{eq:5.7}
\end{align}
Suppose there are $v_1(t)$, $v_2(t)\in {H}^{s+1}_x(\R^+)$, $\norm{v_1}_{{H}^{s+1}_x(\R^+)},\norm{v_2}_{{H}^{s+1}_x(\R^+)}\leq {c}<\frac{1}{c_a}$. Then
\begin{multline}
    \norm{u_1(\cdot,t)-u_2(\cdot,t)}_{\tilde{H}^{s+1}_x(\R^+) } \\= \norm{( v_1-\frac{v_1^2}{2}+\frac{v_1^3}{3}-\frac{v_1^4}{4}+...)-( v_2-\frac{v_2^2}{2}+\frac{v_2^3}{3}-\frac{v_2^4}{4}+...)}_{{H}^{s+1}_x(\R^+)}\\
    =\norm{(v_1-v_2)\left (1-\frac{v_1+v_2}{2}+\frac{v_1^2+v_1v_2+v_2^2}{3}-...\right) }_{{H}^{s+1}_x(\R^+)}\\
    \leq c_a \norm{(v_1-v_2) }_{{H}^{s+1}_x(\R^+)}\left (  1+c+c_ac^2+c_a^2c^3+...\right) \lesssim \norm{(v_1-v_2) }_{{H}^{s+1}_x(\R^+)}. \label{eq:5.8}
\end{multline}
Note that by the classical definitions of the exponential and natural logarithm functions,
\begin{equation}
    u_0(x)= \partial_x \ln(1+v_0)= \partial_x\left(v_0-\frac{v_0^2}{2}+\frac{v^3_0}{3}-... \right). \label{eq:5.9}
\end{equation}
Since $v(\cdot,t)\rightarrow v_0(t)$ in ${H}^{s+1}_x(\R^+) $ as $t\rightarrow 0$, and the map $t\mapsto v(t)$ is $C([0,T];H^{s'}_x(0,\infty))$, we also have $u(\cdot,t)\rightarrow u_0(t)$ in ${\tilde{H}}^{s+1}_x(\R^+) $, and that the map $t\mapsto v(t)$ is 
$C([0,T];\tilde{H}^{s+1}_x(0,\infty))$ (Theorem \ref{thm:2.1}).
\subsection{Time estimate}\label{subsec:5.2}
Combining Theorem \ref{thm:2.1} and Lemma \ref{lem:4.2}, for sufficiently small upper bounds $c_1$ of $\norm{u_0}_{\tilde{H}^s_x}$ and $c_2$ of $\norm{g_0}_{H^{\frac{2s+1}{4}}_t}$, 
\begin{equation}
    \norm{v(x,\cdot)}_{H^{\frac{2s+3}{4}}_t(0,T)} \leq c<\frac{1}{c_a}, \label{eq:5.10}
\end{equation}
and
\begin{equation}
    \norm{v_x(x,\cdot)}_{H^{\frac{2s+1}{4}}_t(0,T)} \leq c<\frac{1}{c_a}, \label{eq:5.11}
\end{equation}
Recalling equation \eqref{eq:5.2}
\begin{equation*}
    u(x,t)= ( v_x-vv_x+v^2v_x-...)=v_x(1-v+v^2-...).
\end{equation*}
The function $(1-v(x)+v(x)^2-...)$ converges to a function $r(x,\cdot) \in H^{\frac{2s+3}{4}}_t(0,T) \subset H^{\frac{2s+1}{4}}_t(0,T)$. $v_x(x,\cdot) \in H^{\frac{2s+1}{4}}_t(0,T)$ Thus, $u(x,\cdot)$ is a well-defined function in $H^{\frac{2s+1}{4}}_t(0,T)$. Moreover, from equation \eqref{eq:5.2}, we have the bound
\begin{equation}
    \norm{u(x,\cdot)}_{H^{\frac{2s+1}{4}}_t(0,T)} \leq cc_a \left( 1+\frac{c}{1-cc_a} \right). \label{eq:5.12}
\end{equation}
To prove that $u(x,\cdot)$ depends continuously on $u_0$ and $g_0$, it suffices to prove that $u(x,\cdot)$ depends continuously on $v(x,\cdot)$ and $v_x(x,\cdot)$.
Let
\begin{align}
    u_1(x,\cdot)&=v_{1,x}(1-v_1+v_1^2-...), \label{eq:5.13}\\
     u_2(x,\cdot)&=v_{2,x}(1-v_2+v_2^2-...). \label{eq:5.14}
\end{align}
Then
\begin{multline}
    \norm{u_1(x,\cdot)-u_2(x,\cdot)}_{H^{\frac{2s+1}{4}}_t(0,T)} \\=\norm{v_{1,x}(1-v_1+v_1^2-...)-v_{2,x}(1-v_2+v_2^2-...)}_{H^{\frac{2s+1}{4}}_t(0,T)} \\
    \leq \norm{v_{1,x}-v_{2,x} }+ \norm{v_{1,x}v_1-v_{2,x}v_2}+\norm{v_{1,x}v_1^2-v_{2,x}v_2^2}+... \\
    \leq \norm{v_{1,x}-v_{2,x} }+(\norm{v_{1,x}(v_1-v_2)+(v_{1,x}-v_{2,x})v_2}+...\\+(\norm{v_{1,x}(v^n_1-v^n_2)+(v_{1,x}-v_{2,x})v^n_2}+...\\
    \leq \norm{v_{1,x}-v_{2,x} }+(cc_a\norm{v_1-v_2 }
    +cc_a\norm{v_{1,x}-v_{2,x}})+...\\+(n(cc_a)^n\norm{v_1-v_2 }
    +(cc_a)^n\norm{v_{1,x}-v_{2,x}})+...\\
    \leq \norm{v_{1,x}-v_{2,x} }(1+cc_a+(cc_a)^2+...)+\norm{v_{1}-v_{2} }(cc_a+2(cc_a)^2+3(cc_a)^3+...)\\
    \lesssim \norm{v_{1,x}(x)-v_{2,x}(x) }_{H^{\frac{2s+1}{4}}_t(0,T)}+\norm{v_{1}(x)-v_{2}(x) }_{H^{\frac{2s+1}{4}}_t(0,T)}. \label{eq:5.15}
\end{multline}
 From the fact that $v(x)$ and $v_x(x)$ are continuous in ${H}^{\frac{2s+3}{4}}_t(0,T) $ and ${H}^{\frac{2s+1}{4}}_t(0,T) $, the map $x\mapsto u(x)$ is $C([0,\infty);{H}^{\frac{2s+1}{4}}_t(0,T))$. As $x\rightarrow 0^+$, $u(x,t)\rightarrow u(0,t)= h_1(t) (1-v_0(t)+v_0^2(t)-...).$ Since $h_1(t)= g_0(t)[1+v(0,t)]$ (in the case of smooth function), $u(x,t)\rightarrow g_0$ in $H^{\frac{2s+1}{4}}_t(0,T)$, thus satisfies the boundary condition

 The argument that $u(x,t)$ is a weak solution of the Burgers' equation follows naturally from smooth estimate of the initial and boundary data.
\section{Analysis of the Initial Data under Condition 2}\label{sec:6}
\begin{prop}\label{prop:6.1}
    The initial data $v_0$ in case 2 is in $H^{s+1}([0,\infty))$, with the lower bound $e^{-c_ac_3} -1 >-1$.
\end{prop}
\begin{proof}
 Applying the same analysis as in section 3, $v_0$ is in $H^{s+1}([0,\infty))$, and depends continuously on $u_0$. Since $1/2<s<3/2$, $v_0$ is $C^1$.
From the upper bound $\norm{u_0}_{ \tilde{H}^{s+1}} \leq c_3$, we get $\norm{\int^x_\infty u_0}_{L^{\infty}} \leq c_ac_3$ (Sobolev embedding). Hence, we deduce the lower bound 
    \begin{equation}
    v_0=e^{\int^x u_0} -1 \geq e^{-c_ac_3} -1 >-1. \label{eq:6.1}
\end{equation}
\end{proof}
\begin{remark}
    The lower bound $e^{-c_ac_3} -1 >-1$ of the initial data $v_0$ is significant because the Cole--Hopf transformation in the inverse direction, from the heat equation back to the Burgers' equation, is
    \begin{equation*}
        u(x,t)=\frac{v_x(x,t)}{1+v(x,t)}.
    \end{equation*}
    If we can establish that $v(x,t)>1$ in $[0,\infty)\times[0,T]$, then the inverse map is well defined everywhere. In Condition 2, we will use the strong minimum principle for the heat equation, to get the lower bound
    \begin{equation}
        \min_{[0,\infty)\times[0,T]} v(x,t) =\min_{[0,\infty)} v_0(x) \geq e^{-c_ac_3} -1 >-1. \label{eq:6.2}
    \end{equation}
    The strong minimum principle that we need is stated below.
\end{remark}
\begin{lemma} \label{lem:6.3}
Let $v(x,t)$ be a classical solution to the heat equation on the half-line with Dirichlet initial condition and Neumann boundary condition
\begin{align}
    v(x,0)&=v_0(x), \label{eq:6.3} \\
    v(0,t)&=h_1(t), \label{eq:6.4}
\end{align}
where $v_0(x)\in C([0,\infty),\R)$, $\lim_{x\rightarrow\infty} v_0(x)=0$; $h_1(t)\in C([0,T],\R)$, and $h_1\leq 0$ in $[0,T]$. 
We require $v(x,t)$ to decay to $0$ as $x\rightarrow \infty$.
Then we have
\begin{equation}
     \min_{[0,\infty)\times[0,T]} v(x,t) =\min_{[0,\infty)} v_0(x). \label{eq:6.5}
\end{equation}
\end{lemma}
\begin{remark}
    The above minimum principle is strong, since it rules out the possibility of $v$ attaining its strict minimum on the boundary $\{ (0,t),t\in (0,T]\} $. 
\end{remark}
\begin{proof}[Proof of Lemma 6.3]
    Consider the parabolic boundary 
    \begin{equation}
        \Sigma= \{ (0,t),\ t\in (0,T]\} \cup \{ (x,0),\ x\in [0,\infty)\} \cup \{ (\infty,t),\ t\in [0,T] \}. \label{eq:6.6}
    \end{equation}
    The weak minimum principle establishes that 
    \begin{equation}
        \min_{[0,\infty)\times[0,T]} v(x,t) =\min_{\Sigma} v(x,t). \label{eq:6.7}
    \end{equation}
    But easily observed that
    \begin{equation}
        \min_{\{ (x,0),\ x\in [0,\infty)\} \cup \{ (\infty,t),\ t\in [0,T] \}} v(x,t) =\min v_0(x), \label{eq:6.8}
    \end{equation}
    hence, if 
    \begin{equation}
        \min_{\{(0,t),\ t\in (0,T]\}} v(x,t) \geq\min v_0(x), \label{eq:6.9}
    \end{equation}
    then we are done.

\noindent Assume that 
\begin{equation}
        \min_{\{(0,t),\ t\in (0,T]\}} v(x,t) <\min v_0(x), \label{eq:6.10}
    \end{equation}
then there must exist a smallest point $0<s\leq T$ such that 
\begin{equation}
    v(0,s)<v(0,t),\quad 0\leq t<s. \label{eq:6.11}
\end{equation}
If $v$ has an interior minimum, then $v$ must be constant (see for example Lemma 2.7 of \cite{Lieberman}). A contradiction. Hence, 
\begin{equation*}
    v(0,s)<v(x,t),\ \text{ for all $(x,t)$ such that } x\in[0,\infty),\ \ 0\leq t<s.
\end{equation*}
Applying Lemma 2.6 of \cite{Lieberman}, which is in a sense the adaptation of the Hopf boundary lemma to parabolic equations, we must have
\begin{equation}
    h_1(s)=\partial_x v(s,0) >0. \label{eq:6.12}
\end{equation}
But this contradicts the assumption that $h_1 \leq 0$ on $[0,T]$.
\end{proof}

\section{Boundary data under Condition 2}\label{sec:7}
From the relationship
\begin{equation*}
    h_1(t)=g_0[1+v(0,t)],
\end{equation*}
we observe the following. If we have $v(x,t)>-1$ in the entire bounded domain, then from the assumption $g_0\leq 0$, $h_1$ would be less than or equal to zero on $(0,T)$. Conversely, if $h_1 \leq 0$, then by Lemma \ref{lem:6.3}, $v(x,t)\geq \min v_0 >-1$. We will integrate this observation into the unique $H^{\frac{2s+1}{4}}_t (0,T)$ solution of the Volterra equation below, to show that this solution must be less than or equal to zero.

\noindent As in section \ref{sec:4}, let 
\begin{equation}
    F(t):=\theta(t) \bigg[g_0(t)+\frac{g_0(t)}{2\pi} \int_\R e^{-k^2t} \hat{V_0}(k) dk - V_x(0,t) \bigg], \label{eq:7.1}
\end{equation}
\begin{equation}
    {h^\ddagger}(t)=(h_1(t)-V_x(0,t))\theta(t), \label{eq:7.2}
\end{equation}
and the Volterra equation
\begin{equation}
    {h^\ddagger}(t)
    = F(t)-\frac{g_0(t)\theta(t)}{\pi} \int_0^t  {h^\ddagger}(s)\frac{1}{\sqrt{t-s}}  ds. \label{eq:7.3} 
\end{equation}
Multiplying by $e^{-Mt}$ on both side, for $M$ to be determined, we have
\begin{equation}
    e^{-Mt} {h^\ddagger}(t)
    = e^{-Mt}F(t)-\frac{g_0(t)\theta(t)}{\pi} \int_0^t \bigg (e^{-Ms}{h^\ddagger}(s)\bigg )  \left( e^{-M(t-s)}\frac{1}{\sqrt{t-s}} \right)  ds. \label{eq:7.4} 
\end{equation}

\begin{prop}\label{prop:7.1}
    Given $f(t) \in H^{\frac{2s+1}{4}}_{0,t} (0,2)$, the map
    \begin{equation}
        R: f(t)\mapsto-\frac{g_0(t)\theta(t)}{\pi} \int_0^t  f(s) \left( e^{-M(t-s)}\frac{1}{\sqrt{t-s}} \right) ds , \label{eq:7.5}
    \end{equation}
 maps $f$ to the function $Rf \in H^{\frac{2s'-1}{4}}_{0,t} (0,2)$. 
    Furthermore,
    \begin{equation}
       \norm{Rf}_{H^{\frac{2s+1}{4}}_{t} (0,2)} \lesssim \norm{e^{-Mt}\frac{1}{\sqrt{t}} }_{L^1} \norm{g_0}_{H^{\frac{2s+1}{4}}} \norm{f}_{H^{\frac{2s+1}{4}}_{t} (0,2)}. \label{eq:7.6} 
    \end{equation}
Thus, there is a sufficiently large $M$ such that $R$ is a contraction map of $f$.
\end{prop}
\begin{proof}
    Similar to Proposition \ref{prop:4.3}, but with the $L^1$ function $ \frac{1}{\sqrt{t}}$ replaced by $e^{-Mt}\frac{1}{\sqrt{t}} $.
\end{proof}
\begin{lemma} \label{lem:7.2}
    The Volterra equation \eqref{eq:7.3} has a unique $h^\ddagger$ solution in ${H^{\frac{2s+1}{4}}_{t} (0,2)}$. Furthermore, this solution satisfy
    \begin{equation}
        h^\ddagger(t)+V_x(t) \leq 0,\quad t\in [0,T], \label{eq:7.7}
    \end{equation}
    in addition to depend continuously on $u_0$ and $g_0$.
\end{lemma}
\begin{proof}
Equation \eqref{eq:7.2} implies back 
\begin{align}
    h^{\ddagger}(t)+V(0,t)_x&= g_0(t)\left(1+S[0,h^{\ddagger}](0,t)+S[V_0](0,t)\right) \nonumber \\
    &= g_0(t)\left (1+S[v_0, h^{\ddagger}+V_x](0,t) \right), \label{eq:7.8}
\end{align}
for $x\geq 0$, $t\in (0,T)$ (since $\theta=1$ on $[0,T]$). The proof follows closely that of Lemma \ref{lem:4.2} and Proposition \ref{prop:4.5}.

To prove inequality \eqref{eq:7.7}, suppose we start the iteration process with $h^{\ddagger}_{(0)}$ such that 
\begin{equation}
   h_{1,(0)}= h^{\ddagger}_{(0)}+V_x(0,t) \leq 0\quad\text{on $[0,T]$}. \label{eq:7.9}
\end{equation}
From Lemma \ref{lem:2.4}, $S[v_0, h^{\ddagger}_{(0)}+V_x](x,t)$ is a classical solution of the Neumann BC heat equation on the half-line. From equation (7.9), we can apply the strong minimum principle (Lemma 6.3) in $(0,\infty)\times [0,T]$.
Thus, 
\begin{equation*}
    1+ S[v_0, h^{\ddagger}_{(0)}+V_x](0,t) \geq 1+ \min_{x} v_0\geq e^{-c_ac_3},\quad t\in [0,T].
\end{equation*}
Since in Condition 2, $g_0(t) \leq 0 $ for $ t\in [0,2]$, we deduce that the RHS of (7.8) is less than or equal to zero on $[0,T]$ as well, which implies
\begin{equation*}
   h_{1,(1)}= h^{\ddagger}_{(1)}+V_x(0,t) \leq 0\quad\text{on $[0,T]$}.
\end{equation*}
By induction, this holds true for all iteration terms.
Convergence in $H^{\frac{2s+1}{4}}(0,2) $, for $\frac{2s+1}{4}>\frac{1}{2}$, implies uniform convergence in $[0,T]$. We easily deduce that the unique solutions $h^\ddagger(t) \in  H^{\frac{2s+1}{4}}_0(0,2)$, and the boundary condition $h_1(t)$, satisfy
\begin{equation*}
   h_{1}= h^{\ddagger}_{}+V_x(0,t) \leq 0\quad\text{on $[0,T]$},
\end{equation*}
in addition to be continuously dependent on $u_0$ and $g_0$, as in section 4.
\end{proof}
\section{well--posedness under Condition 2}\label{sec:8}
From Lemma \ref{lem:7.2}, we have $h_{1}(t)\leq 0$ for $t\in[0,T]$. From Lemma \ref{lem:2.4}, the
solution $v(x,t)=S[v_{0},h_{1}]$ is a classical solution in
$(0,\infty)\times(0,T]$, and applying Lemma \ref{lem:6.3} (the strong minimum principle),
we get that
\begin{equation}
v(x,t)\geq \min_{x\in(0,\infty)}v_{0}(x)\geq e^{-c_{a}c_{3}}-1>-1, \label{eq:8.1}
\end{equation}
for all $(x,t)\in[0,\infty)\times[0,T]$. We now examine the return map
\[
u(x,t)=\partial_x\ln (v+1)=\frac{v_{x}}{v+1}.
\]
\subsection{Space estimate}\label{subsubsec:8.0.1}
To establish that $u(x,t)\in \tilde{H}^{s+1}(\R^+)$, we need to show that 
\begin{equation*}
    \norm{\ln(v+1)}_{H^{s+1}(\R^+)} = \norm{\ln(v+1)}_{L^2(\R^+)}+\norm{\frac{v_x}{v+1}}_{H^s(\R^+)}.
\end{equation*}

For $0<\beta<1$, the Sobolev seminorm $[\cdot]_\beta$ is defined as
\[
[w]_{\beta}^{2}
:=\int_{0}^{\infty}\!\!\int_{0}^{\infty}
\frac{|w(x)-w(y)|^{2}}{|x-y|^{1+2\beta}}\,dx\,dy ,
\]
so that $\|w\|_{H^{\beta}(\mathbb{R}^{+})}^{2}\simeq
\|w\|_{L^{2}(\mathbb{R}^{+})}^{2}+[w]_{\beta}^{2}$ for $0<\beta<1$,
$\|w\|_{H^{1}(\mathbb{R}^{+})}^{2}=\|w\|_{L^{2}}^{2}+\|w_{x}\|_{L^{2}}^{2}$, and
$\|w\|_{H^{s}(\mathbb{R}^{+})}^{2}\simeq
\|w\|_{H^{1}(\mathbb{R}^{+})}^{2}+[w_{x}]_{s-1}^{2}$ for $1<s<\tfrac32$.
\begin{prop}\label{prop:8.1}
    Let $a$, $b$ be two positive real number such that $a,b\geq e^{-c_ac_3} >0$. Then
    \begin{equation}
        |\ln(a)-\ln(b)|\leq e^{c_ac_3} |a-b|. \label{eq:8.2}
    \end{equation}
\end{prop}
\begin{proof}
    Without loss of generality, let $a\leq b$. By the mean value theorem, there is a $c\in (a,b)$ such that $|\ln(a)-\ln(b)|\leq \frac{1}{c} |a-b|$. Since $c\in (a,b)$, $\frac{1}{c} \leq e^{c_ac_3}$.
\end{proof}
For any $x\in \R^+$, let $a=\ln(v(x)+1)$, $b=\ln(1)$, from Proposition \ref{prop:8.1}, we have
\begin{equation}
    |\ln(v(x)+1)|=|\ln(v(x)+1)-\ln(1)| \leq e^{c_ac_3} |v(x)|. \label{eq:8.3}
\end{equation}
Thus, 
\begin{equation}
    \norm{\ln(v+1)}_{L^2(\R^+)}\leq e^{c_ac_3} \norm{v(x)}_{L^2(\R^+)} \leq e^{c_ac_3} \norm{v(x)}_{H^{s+1}(\R^+)} . \label{eq:8.4}
\end{equation}
To estimate $\norm{\frac{v_x}{v+1}}_{H^s(\R^+)}$, we need the below lemma.
\begin{lemma}[Quotient lemma]\label{lem:quotient}
Let $\tfrac12<s<\tfrac32$ and let $f,g\in H^{s}(\mathbb{R}^{+})$ be real valued.
Assume that
\[
G:=1+g\geq c>0 \qquad \text{on } \mathbb{R}^{+} .
\]
Then
\[
u:=\frac{f}{1+g}\in H^{s}(\mathbb{R}^{+}),
\]
and
\begin{equation}\label{eq:quotient}
\|u\|_{H^{s}(\mathbb{R}^{+})}
\leq C(s,c)\,\|f\|_{H^{s}(\mathbb{R}^{+})}
\bigl(1+\|g\|_{H^{s}(\mathbb{R}^{+})}\bigr)^{3},
\end{equation}
where $C(s,c)>0$ depends on $s$ and $c$.
\end{lemma}

\begin{proof}
See Appendix \ref{app:B}.
\end{proof}

Take $f=v_{x}(t)$ and $g=v(t)$, then
$1+g=1+v(t)\geq e^{-c_{a}c_{3}}>0$. Since $v(t)\in H^{s+1}(\mathbb{R}^{+})$, that gives
$v_{x}(t)\in H^{s}(\mathbb{R}^{+})$. Lemma \ref{lem:quotient} yields directly
\begin{equation}
    \Bigl\|\frac{v_{x}(t)}{1+v(t)}\Bigr\|_{H^{s}(\mathbb{R}^{+})}
\lesssim C\bigl(e^{-c_{a}c_{3}}\bigr)\,
\|v(t)\|_{H^{s+1}(\mathbb{R}^{+})}
\bigl(1+\|v(t)\|_{H^{s+1}(\mathbb{R}^{+})}\bigr)^{3}. \label{eq:8.6}
\end{equation}
Combining equations \eqref{eq:8.4} and \eqref{eq:8.6}, 
\begin{equation}
    \norm{\ln(v(t)+1)}_{H^{s+1}(\R^+)} \lesssim C\bigl(e^{-c_{a}c_{3}}\bigr)\,
\|v(t)\|_{H^{s+1}(\mathbb{R}^{+})}
\bigl(1+\|v(t)\|_{H^{s+1}(\mathbb{R}^{+})}\bigr)^{3}. \label{eq:8.7}
\end{equation}

Furthermore, for $t_1,t_2\in [0,T]$,
\begin{equation}
    \norm{u(t_1)-u(t_2)}_{\tilde{H}^{s+1}}\simeq \norm{\ln(v(t_1)+1)-\ln(v(t_2)+1)}_{L_2}+ \norm{\frac{v_x(t_1)}{v(t_1)+1}- \frac{v_x(t_2)}{v(t_2)+1}}_{H^s}. \label{eq:8.8}
\end{equation}
    Similar to equation \eqref{eq:8.4}, 
    \begin{equation}
        \norm{\ln(v(t_1)+1)-\ln(v(t_2)+1)}_{L_2}\leq e^{c_ac_3} \norm{v(t_1)-v(t_2)}_{L^2}\lesssim e^{c_ac_3} \norm{v(t_1)-v(t_2)}_{H^{s+1}}. \label{eq:8.9}
    \end{equation}
Also,
    \begin{multline}
    \norm{\frac{v_x(t_1)}{v(t_1)+1}- \frac{v_x(t_2)}{v(t_2)+1}}_{H^s}=\norm{\frac{v_x(t_1)(v(t_2)-v(t_1))-(v_x(t_2)-v_x(t_1))(v(t_1)+1)}{ (v(t_2)+1)(v(t_1)+1)} } _{H^s} \\
    \lesssim \norm{\frac{v_x(t_1)(v(t_2)-v(t_1))}{ (v(t_2)+1)(v(t_1)+1)} } _{H^s}+\norm{\frac{(v_x(t_2)-v_x(t_1))(v(t_1)+1)}{ (v(t_2)+1)(v(t_1)+1)} } _{H^s} \\
    \lesssim \norm{\frac{v_x(t_1)}{v(t_1)+1} }_{H^s} \norm{\frac{v(t_2)-v(t_1)}{v(t_2)+1}}_{H^{s}} + \norm{\frac{v_x(t_2) -v_x(t_1)}{v(t_2)+1}}_{H^s}. \label{eq:8.10}
\end{multline}
By equation \eqref{eq:8.6}, $\norm{\frac{v_x(t_1)}{v(t_1)+1} }_{H^s}$ is bounded, and by Lemma \ref{lem:quotient}, as $t_2\rightarrow t_1$, $\norm{\frac{v(t_2)-v(t_1)}{v(t_2)+1}}_{H^{s}}$, $\norm{\frac{v_x(t_2) -v_x(t_1)}{v(t_2)+1}}_{H^s}$ goes to zero, hence $\norm{\frac{v_x(t_1)}{v(t_1)+1}- \frac{v_x(t_2)}{v(t_2)+1}}_{H^s}$ goes to zero ($v$ is continuous in $H^{s+1}$ with respect to $t$). Similarly, $\norm{\ln(v(t_1)+1)-\ln(v(t_2)+1)}_{L_2}$ would also goes to zero. Combining together, $\norm{u(t_1)-u(t_2)}_{\tilde{H}^{s+1}}$ must also go to zero. Hence, $u(t)$ is $C([0,T],\tilde{H}^{s+1})$. In particular, $u(t)\rightarrow u_0(x)$ in $\tilde{H}^{s+1}(\R^+)$ as $t\rightarrow 0$.
\subsection{Time estimate}\label{subsubsec:8.0.2}
Restricting the domain to $(0,T]$ and the regularity range to $(\frac{1}{2},1)$, similar to case 1 of Lemma \ref{lem:quotient}, we have
\begin{lemma} \label{lem:8.3}
    Let $\frac12<b<1$ and let $f,g\in H^{b}((0,T])$ be real valued.
Assume that
\[
G:=1+g\geq c>0 \qquad \text{on } (0,T] .
\]
Then
\[
u:=\frac{f}{1+g}\in H^{b}((0,T]),
\]
and
\begin{equation}\label{eq:quotient}
\|u\|_{H^{b}((0,T])}
\;\lesssim\; c^{-2}\,\|f\|_{H^{b}((0,T])}
\bigl(1+\|g\|_{H^{b}((0,T])}\bigr).
\end{equation}
\end{lemma}

Since $\frac{1}{2}<s<\frac{3}{2}$, $\frac{1}{2}<\frac{2s+1}{4}<1$.
For any $x\in [0,\infty)$,
\begin{multline}
    \norm{u(x)}_{H^{\frac{2s+1}{4}}_t} = \norm{\frac{v_x(x)}{v(x)+1}}_{H^{\frac{2s+1}{4}}_t}  \lesssim e^{2c_ac_3} \norm{{v_x(x)}}_{H^{\frac{2s+1}{4}}_t} \left(1+\norm{ v(x)}_{ H^{\frac{2s+1}{4}}_t}\right) \\
    \lesssim e^{2c_ac_3} \norm{{v_x(x)}}_{H^{\frac{2s+1}{4}}_t} \left(1+\norm{ v(x)}_{ H^{\frac{2s+3}{4}}_t}\right). \label{eq:8.12}
\end{multline}
For $x_1$, $x_2 \in [0,\infty)$,
\begin{multline}
    \norm{u(x_1)-u(x_2)}_{H^{\frac{2s+1}{4}}_t}=\norm{\frac{v_x(x_1)}{v(x_1)+1}-\frac{v_x(x_2)}{v(x_2)+1}}_{H^{\frac{2s+1}{4}}_t}\\
    =\norm{\frac{v_x(x_1)(v(x_2)-v(x_1))-(v_x(x_2)-v_x(x_1))(v(x_1)+1)}{ (v(x_2)+1)(v(x_1)+1)} } _{H^{\frac{2s+1}{4}}_t} \\
    \lesssim \norm{\frac{v_x(x_1)(v(x_2)-v(x_1))}{ (v(x_2)+1)(v(x_1)+1)} } _{H^{\frac{2s+1}{4}}_t}+\norm{\frac{(v_x(x_2)-v_x(x_1))(v(x_1)+1)}{ (v(x_2)+1)(v(x_1)+1)} } _{H^{\frac{2s+1}{4}}_t} \\
    \lesssim \norm{\frac{v_x(x_1)}{v(x_1)+1} }_{H^{\frac{2s+1}{4}}_t} \norm{\frac{v(x_2)-v(x_1)}{v(x_2)+1}}_{H^{\frac{2s+1}{4}}_t} + \norm{\frac{v_x(x_2) -v_x(x_1)}{v(x_2)+1}}_{H^{\frac{2s+1}{4}}_t}. \label{eq:8.13}
\end{multline}
Since $v(x)$ and $v_x(x)$ are continuous in $H^{\frac{2s+3}{4}}_t$ and $H^{\frac{2s+1}{4}}_t$, applying Lemma \ref{lem:8.3} thrice, as $x_2\rightarrow x_1$, $\norm{u(x_1)-u(x_2)}_{H^{\frac{2s+1}{4}}_t}$ goes to zero. Hence, we have $u(x)$ is $C([0,\infty), H^{\frac{2s+1}{4}}_t$). In particular, $u(x)\rightarrow g_0$ as $x\rightarrow0$.
\appendix
\section{Proof of Linear Estimates}\label{sec:A}
We cite the following lemmas from \cite{heat}.

\begin{lemma}[\cite{heat}, Lemma 2.1]
If $\gamma \in \mathbb{C}$ with $\operatorname{Im}(\gamma) > 0$, then
\[
|e^{i\gamma k x} - e^{i\gamma k z}| \leq \sqrt{2} \left(1 + \frac{|\operatorname{Re}(\gamma)|}{\operatorname{Im}(\gamma)} \right) |e^{-\operatorname{Im}(\gamma) k x} - e^{-\operatorname{Im}(\gamma) k z}| \quad \forall k,x,z \geq 0.
\]
\end{lemma}

\begin{lemma}[$L^2$-boundedness of the Laplace transform]
The map
\[
\mathcal{L}: \phi \mapsto \int_0^\infty e^{-\tau t} \phi(\tau) \, d\tau
\]
is bounded from $L^2_\tau(0,\infty)$ into $L^2_t(0,\infty)$ with
\[
\|\mathcal{L}\{\phi\}\|_{L^2_t(0,\infty)} \leq \sqrt{\pi} \, \|\phi\|_{L^2_\tau(0,\infty)}.
\]
\end{lemma}
We also need the following equation for $k\in \R$:
\begin{equation}
\int_0^\infty \frac{(1 - e^{-k z})^2}{z^{1+2\beta}} \, dz = k^{2\beta} \int_0^\infty \frac{(1 - e^{-\zeta})^2}{\zeta^{1+2\beta}} \, d\zeta \simeq k^{2\beta}, \quad 0 < \beta < 1. \label{eq:2.9}
\end{equation}
\subsection{Proof of Theorem 2.2 (Linear IVP Estimates)}\label{subsec:A.1}
\textit{Proof.} The space estimate (2.14) holds for all $s \in \mathbb{R}$ (see \cite{heat}).

For the time estimates (2.15) and (2.16), let $m=\frac{2s'+1}{4}$. Since $\frac{3}{2} < s' < \frac{5}{2}$, we have $1 < m < \frac{3}{2}$. let $m_x = \frac{2s'-1}{4}$. we have $\frac{1}{2} < m_x < 1$.

We employ the physical space definition of the Sobolev norm:
\begin{equation}
\|V(x)\|_{H^m_t(0,T)} = \|\partial_t V(x)\|_{L^2_t(0,T)} +\|\partial_t V(x)\|_{L^2_t(0,T)} + \|\partial_t V(x)\|_{m-1}, \label{eq:A.2} 
\end{equation}
where the fractional norm is defined as
\begin{equation*}
\|V(x)\|_{(m-1)} = \left( \int_{t=0}^T \int_{z=0}^{T-t} \frac{|V(x,z+t) - V(x,t)|^2}{z^{1+2(m-1)}} \, dz \, dt \right)^{\frac{1}{2}}. 
\end{equation*}
Similarly,
\begin{equation}
\|V_x(x)\|_{H^{m_x}_t(0,T)} =\| V_x(x)\|_{L^2_t(0,T)} + \| U_x(x)\|_{m_x}, \label{eq:A.3} 
\end{equation}
where the fractional norm is defined as
\begin{equation*}
\|V_x(x)\|_{m_x} = \left( \int_{t=0}^T \int_{z=0}^{T-t} \frac{|V_x(x,z+t) - V_x(x,t)|^2}{z^{1+2{m_x}}} \, dz \, dt \right)^{\frac{1}{2}}. 
\end{equation*}
Consider the formula of $V$, $V_t$, $V_x$ by differentiating (2.11),
\begin{align}
\|V(x)\|_{L^2_t(0,T)}^2 &\approx \int_0^T \left[ \int_0^\infty  e^{-\xi^2 t} (|\widehat{V}_0(-\xi)| + |\widehat{V}_0(\xi)|) \, d\xi \right]^2 dt . \label{eq:A.4}\\
\|V_t(x)\|_{L^2_t(0,T)}^2 &\approx \int_0^T \left[ \int_0^\infty  e^{-\xi^2 t} \xi^2 (|\widehat{V}_0(-\xi)| + |\widehat{V}_0(\xi)|) \, d\xi \right]^2 dt . \label{eq:A.5}\\
\|V_x(x)\|_{L^2_t(0,T)}^2 &\approx \int_0^T \left[ \int_0^\infty  e^{-\xi^2 t} \xi (|\widehat{V}_0(-\xi)| + |\widehat{V}_0(\xi)|) \, d\xi \right]^2 dt. \label{eq:A.6}
\end{align}
For $j=0,1,2$, we have
\begin{subequations} \label{eq:A.7}
    \begin{align}
        \int_0^T \left[ \int_0^\infty  e^{-\xi^2 t} \xi^j (|\widehat{V}_0(-\xi)| + |\widehat{V}_0(\xi)|) \, d\xi \right]^2 dt \label{eq:A.7a} \\ \lesssim \int_0^T \left[ \int_0^1\xi^j e^{-\xi^2 t} (|\widehat{V}_0(-\xi)| + |\widehat{V}_0(\xi)|) \, d\xi \right]^2 dt \label{eq:A.7b} \\
 + \int_0^T \left[ \int_1^\infty \xi^j e^{-\xi^2 t} (|\widehat{V}_0(-\xi)| + |\widehat{V}_0(\xi)|) \, d\xi \right]^2 dt. \label{eq:A.7c} 
    \end{align}
\end{subequations}

By the Cauchy-Schwarz inequality, we find
\begin{equation}
(A.7b) \lesssim \int_0^T \left(  \int_0^1 e^{-2\xi^2 t} (1 + \xi^2)^{-s'}\xi^{2j} \, d\xi \right) \|V_0\|_{H^{s'}_x(\mathbb{R})}^2 \, dt \leq c_s T \|V_0\|_{H^{s'}_x(\mathbb{R})}^2, \quad {s'} \in \mathbb{R}. \label{eq:3.7}
\end{equation}
Moreover, making the change of variable $\xi = \sqrt{\tau}$, we have
\begin{equation}
(A.7c) = \int_0^T \left( \int_1^\infty  e^{-\tau t} \frac{|\widehat{V}_0(-\sqrt{\tau})| + |\widehat{V}_0(\sqrt{\tau})|}{2 \sqrt{\tau}} \tau^{j/2} \, d\tau \right)^2 dt \lesssim \|\mathcal{L}\{\phi\}\|_{L^2_t(0,\infty)}^2, \label{eq:3.8}
\end{equation}
where $\mathcal{L}\{\phi\}(t) = \int_0^\infty e^{-\tau t} \phi(\tau) \, d\tau$ is the Laplace transform of the function
\[
\phi(\tau) = 
\begin{cases}
 \tau^{\frac{j-1}{2}} (|\widehat{V}_0(-\sqrt{\tau})| + |\widehat{V}_0(\sqrt{\tau})|), & \tau \geq 1, \\
0, & 0 \leq \tau < 1.
\end{cases}
\]
Hence, Lemma 2.2 yields
\[
(A.7c) \lesssim \int_1^\infty \tau^{j-1} |\widehat{V}_0(-\sqrt{\tau})|^2 \, d\tau + \int_1^\infty \tau^{j-1}|\widehat{V}_0(\sqrt{\tau})|^2 \, d\tau,
\]
and, letting $\xi = -\sqrt{\tau}$ and $\xi = \sqrt{\tau}$ in the first and the second integral, respectively, we obtain
\begin{equation}
(A.7c) \lesssim \int_{|\xi| \geq 1} |\xi|^{2j-1} |\widehat{V}_0(\xi)|^2 \, d\xi \leq \|V_0\|_{H^{3/2}_x(\mathbb{R})}^2 \leq \|V_0\|_{H^s_x(\mathbb{R})}^2, \quad {s'} \geq \frac{3}{2}. \label{eq:3.9}
\end{equation}
Combining estimates~\eqref{eq:3.7} and~\eqref{eq:3.9}, we find
\begin{align}
    \|V(x)\|_{L^2_t(0,T)} &\lesssim \|V_0\|_{H^{s'}_x(\mathbb{R})}, \quad {s'} \geq -\frac{1}{2}, \; x \in \mathbb{R}. \label{eq:3.10} \\
    \|V_t(x)\|_{L^2_t(0,T)} &\lesssim \|V_0\|_{H^{s'}_x(\mathbb{R})}, \quad {s'} \geq \frac{3}{2}, \; x \in \mathbb{R}. \label{eq:A.12}  \\
    \|V_x(x)\|_{L^2_t(0,T)} &\lesssim \|V_0\|_{H^{s'}_x(\mathbb{R})}, \quad {s'} \geq \frac{1}{2}, \; x \in \mathbb{R}. \label{eq:A.13} 
\end{align}

\textbf{Estimate of $\|V_x(x)\|_{m_x}$:}
    For the fractional norm $\|V_x(x)\|_{m_x}$, starting from formula (2.12), we compute
\begin{multline}
    |V_x(x,t+z) - V_x(x,t)| \lesssim \int_0^1 e^{-\xi^2 t} (1 - e^{-\xi^2 z}) \xi (|\widehat{V}_0(-\xi)| + |\widehat{V}_0(\xi)|) \, d\xi \\+ \int_1^\infty e^{-\xi^2 t} (1 - e^{-\xi^2 z}) \xi (|\widehat{V}_0(-\xi)| + |\widehat{V}_0(\xi)|) \, d\xi. \label{eq:A.14}
\end{multline}
Hence, from the definition of $\|V(x)\|_{m_x}$, we have
\begin{align}
\|V_x(x)\|_{m_x}^2 &\lesssim \int_0^T \int_0^{T-t} \frac{1}{z^{1+2{m_x}}} \left[ \int_0^1 e^{-\xi^2 t} (1 - e^{-\xi^2 z}) \xi(|\widehat{V}_0(-\xi)| + |\widehat{V}_0(\xi)|) \, d\xi \right]^2 dz \, dt \label{eq:3.11a} \\
&\quad + \int_0^T \int_0^{T-t} \frac{1}{z^{1+2{m_x}}} \left[ \int_1^\infty e^{-\xi^2 t} (1 - e^{-\xi^2 z}) \xi (|\widehat{V}_0(-\xi)| + |\widehat{V}_0(\xi)|) \, d\xi \right]^2 dz \, dt. \label{eq:3.11b}
\end{align}
For any $s \in \mathbb{R}$, the Cauchy-Schwarz inequality implies
\[
(A.15) \lesssim T \|V_0\|_{H^{s'}_x(\mathbb{R})}^2 \int_0^1 (1 + \xi^2)^{-s'} \xi^2 \int_0^T \frac{(1 - e^{-\xi^2 z})^2}{z^{1+2{m_x}}} \, dz \, d\xi.
\]
Therefore, employing~\eqref{eq:2.9} for $k = \xi^2$, we find
\begin{equation}
(A.15) \lesssim T \|V_0\|_{H^{s'}_x(\mathbb{R})}^2 \int_0^1 (1 + \xi^2)^{-s'} \xi^{4{m_x}+2} \, d\xi = c_s T \|V_0\|_{H^{s'}_x(\mathbb{R})}^2, \quad 
\frac{3}{2} < {s'} < \frac{5}{2}. \label{eq:3.12}
\end{equation}
Furthermore, letting $\xi = \sqrt{\tau}$, we have
\[
(A.16) \lesssim \int_0^T \int_0^{T-t} \frac{1}{z^{1+2{m_x}}} \|\mathcal{L}\{\phi_z\}(t)\|_{L^2_t(0,\infty)}^2 \, dz \, dt,
\]
where
\[
\phi_z(\tau) = 
\begin{cases}
 (1 - e^{-\tau z}) (|\widehat{V}_0(-\sqrt{\tau})| + |\widehat{V}_0(\sqrt{\tau})|), & \tau \geq 1, \\
0, & 0 \leq \tau < 1.
\end{cases}
\]
Thus, by the Laplace transform bound of Lemma 2.2, we find
\[
(A.16) \lesssim \int_1^\infty \xi(|\widehat{V}_0(-\xi)| + |\widehat{V}_0(\xi)|)^2 \int_0^T \frac{(1 - e^{-\xi^2 z})^2}{z^{1+2{m_x}}} \, dz \, d\xi,
\]
and, using~\eqref{eq:2.9} as before, we obtain
\begin{equation}
(A.16) \lesssim \int_1^\infty \xi^{4{m_x} + 1} (|\widehat{V}_0(-\xi)| + |\widehat{V}_0(\xi)|)^2 \, d\xi \lesssim \|V_0\|_{H^{s'}_x(\mathbb{R})}^2, \quad \frac{3}{2} < {s'} < \frac{5}{2}. \label{eq:3.13}
\end{equation}
Combining estimates (A.17) and (A.18), we deduce
\begin{equation}
\|V_x(x)\|_{m_x} \lesssim \|V_0\|_{H^s_x(\mathbb{R})}, \quad \frac{3}{2} < {s'}< \frac{5}{2}, \; x \in \mathbb{R}, \label{eq:3.14}
\end{equation}
which together with estimate (A.13) and the definition (A.3) yield estimate (2.16).

\textbf{Estimate of $\|V_t(x)\|_{m-1}$:}

    For the fractional norm $\|V_t(x)\|_{m-1}$, starting from formula (2.12), we compute
\begin{multline}
    |V_t(x,t+z) - V_t(x,t)| \lesssim \int_0^1 e^{-\xi^2 t} (1 - e^{-\xi^2 z}) \xi^2 (|\widehat{V}_0(-\xi)| + |\widehat{V}_0(\xi)|) \, d\xi \\+ \int_1^\infty e^{-\xi^2 t} (1 - e^{-\xi^2 z}) \xi^2 (|\widehat{V}_0(-\xi)| + |\widehat{V}_0(\xi)|) \, d\xi. \label{eq:A.20}
\end{multline}
Hence, from the definition of $\|V_t(x)\|_{m-1}$, we have
\begin{align}
\|V_t(x)\|_{m-1}^2 &\lesssim \int_0^T \int_0^{T-t} \frac{1}{z^{1+2{(m-1)}}} \left[ \int_0^1 e^{-\xi^2 t} (1 - e^{-\xi^2 z}) \xi^2(|\widehat{V}_0(-\xi)| + |\widehat{V}_0(\xi)|) \, d\xi \right]^2 dz \, dt \label{eq:3.11a} \\
&\quad + \int_0^T \int_0^{T-t} \frac{1}{z^{1+2{(m-1)}}} \left[ \int_1^\infty e^{-\xi^2 t} (1 - e^{-\xi^2 z}) \xi^2 (|\widehat{V}_0(-\xi)| + |\widehat{V}_0(\xi)|) \, d\xi \right]^2 dz \, dt. \label{eq:3.11b}
\end{align}
For any $s' \in \mathbb{R}$, the Cauchy-Schwarz inequality implies
\[
\eqref{eq:3.11a} \lesssim T \|V_0\|_{H^{s'}_x(\mathbb{R})}^2 \int_0^1 (1 + \xi^2)^{-s'} \xi^4 \int_0^T \frac{(1 - e^{-\xi^2 z})^2}{z^{1+2{(m-1)}}} \, dz \, d\xi.
\]
Therefore, employing~\eqref{eq:2.9} for $k = \xi^2$ and $\beta = {m-1}$, we find
\begin{equation}
\eqref{eq:3.11a} \lesssim T \|V_0\|_{H^{s'}_x(\mathbb{R})}^2 \int_0^1 (1 + \xi^2)^{-s'} \xi^{4{m}} \, d\xi = c_{s'} T \|V_0\|_{H^{s'}_x(\mathbb{R})}^2, \quad 
\frac{3}{2} < {s'}< \frac{5}{2}. \label{eq:3.12}
\end{equation}
Furthermore, letting $\xi = \sqrt{\tau}$, we have
\[
\eqref{eq:3.11b} \lesssim \int_0^T \int_0^{T-t} \frac{1}{z^{1+2{(m-1)}}} \|\mathcal{L}\{\phi_z\}(t)\|_{L^2_t(0,\infty)}^2 \, dz \, dt,
\]
where
\[
\phi_z(\tau) = 
\begin{cases}
 \tau^{\frac{3}{2}}(1 - e^{-\tau z}) (|\widehat{V}_0(-\sqrt{\tau})| + |\widehat{V}_0(\sqrt{\tau})|), & \tau \geq 1, \\
0, & 0 \leq \tau < 1.
\end{cases}
\]
Thus, by the Laplace transform bound of Lemma 2.2, we find
\[
\eqref{eq:3.11b} \lesssim \int_1^\infty \xi^3(|\widehat{V}_0(-\xi)| + |\widehat{V}_0(\xi)|)^2 \int_0^T \frac{(1 - e^{-\xi^2 z})^2}{z^{1+2{(m-1)}}} \, dz \, d\xi,
\]
and, using~\eqref{eq:2.9} as before, we obtain
\begin{equation}
\eqref{eq:3.11b} \lesssim \int_1^\infty \xi^{4{(m-1)} + 3} (|\widehat{V}_0(-\xi)| + |\widehat{V}_0(\xi)|)^2 \, d\xi \lesssim \|V_0\|_{H^{s'}_x(\mathbb{R})}^2, \quad \frac{3}{2} < {s'} < \frac{5}{2}. \label{eq:3.13}
\end{equation}

Combining estimates~\eqref{eq:3.12} and~\eqref{eq:3.13}, we deduce
\begin{equation}
\|V_t(x)\|_{m-1} \lesssim \|V_0\|_{H^{s'}_x(\mathbb{R})}, \quad \frac{3}{2} < {s'} < \frac{5}{2}, \; x \in \mathbb{R}, \label{eq:3.14}
\end{equation}
which together with estimates (A.11), (A.12) and the definition (A.2) yield estimate (2.15).

\subsection{Proof of Theorem 2.3 (Estimates for the Boundary Value Problem Components)}\label{subsec:A.2}

In order to prove the space estimate (2.17), we write $v = v_1 + v_2$ with
\begin{align}
v_1(x, t) &= \int_0^\infty e^{i\gamma_1 k x} H_1(k, t) \, dk, \quad H_1(k, t) =- \frac{1}{\pi} e^{i k^2 t}  \widehat{{h^\dagger}}(k^2), \quad \gamma_1 = a^3 = e^{i 3\pi/4}, \label{eq:A.26} \\
v_2(x, t) &= \int_0^\infty e^{i\gamma_2 k x} H_2(k, t) \, dk, \quad H_2(k, t) = -\frac{1}{\pi} e^{-i k^2 t} \widehat{{h^\dagger}}(-k^2), \quad \gamma_2 = a = e^{i \pi/4}. \label{eq:A.27}
\end{align}
The estimation of $v_1$ and $v_2$ is entirely analogous. Thus, we only provide the details for $v_1$. It is convenient to employ the physical space definition of the $H^{s'}_x(0,\infty)$-norm,
\begin{equation}
\|v_1(t)\|_{H^{s'}_x(0,\infty)} = \sum_{j=0}^{\lfloor {s'} \rfloor} \|\partial_x^j v_1(t)\|_{L^2_x(0,\infty)} + \|\partial_x^{\lfloor {s'} \rfloor} v_1(t)\|_\beta, \quad {s'} = \lfloor {s'} \rfloor + \beta \geq 0, \; 0 \leq \beta < 1, \label{eq:2.7}
\end{equation}
where $\lfloor \cdot \rfloor$ denotes the floor function and, for $0 < \beta < 1$, the fractional norm $\|\cdot\|_\beta$ is defined by
\[
\|v_1(t)\|_\beta^2 = \int_0^\infty \int_0^\infty \frac{|v_1(x,t) - v_1(y,t)|^2}{|x - y|^{1 + 2\beta}} \, dy \, dx \simeq \int_0^\infty \int_0^\infty \frac{|v_1(x+z,t) - v_1(x,t)|^2}{z^{1 + 2\beta}} \, dz \, dx.
\]
There are three cases to consider: (i) $\lfloor {s'} \rfloor = 0$ and $\beta \ne 0$; (ii) $\beta = 0$; (iii) $\lfloor {s'} \rfloor \ne 0$ and $\beta \ne 0$. [Note, the first case is obsolete, since $\frac{3}{2}<s'<\frac{5}{2}$, but it will be more convenient to preserve the structure of the original proof. This actually shows that the relevant estimate is robust for a variety of problems, in various regularity range.]

\textbf{(i) The case $\lfloor {s'} \rfloor = 0$ and $\beta \ne 0$.} Then ${s'} = \beta \in (0,1)$ and we need to estimate $\|v_1(t)\|_{L^2_x(0,\infty)}$ and $\|v_1(t)\|_\beta$. The first norm will be estimated together with the $L^2_x(0,\infty)$-norms of the higher derivatives of $v_1$ in case (ii). For the second norm, we have
\[
\|v_1(t)\|_\beta^2 \leq \int_0^\infty \int_0^\infty \frac{1}{z^{1+2\beta}} \left( \int_0^\infty |e^{i\gamma_1 k(x+z)} - e^{i\gamma_1 k x}| \, |H_1(k,t)| \, dk \right)^2 dx \, dz.
\]
Using Lemma 2.1 with $\gamma = \gamma_1 = e^{i 3\pi/4}$ and making the change of variables $\frac{\sqrt{2}}{2} x \to x$, $\frac{\sqrt{2}}{2} z \to z$, we obtain
\begin{equation}
\|v_1(t)\|_\beta^2 \lesssim \int_0^\infty \frac{1}{z^{1+2\beta}} \int_0^\infty \left( \int_0^\infty e^{-k x} (1 - e^{-k z}) |H_1(k,t)| \, dk \right)^2 dx \, dz. \label{eq:2.8}
\end{equation}
Using Lemma 2.2 for $\phi= (1-e^{-kz)}|H_1(k,t)|$
\begin{equation*}
    \|v_1(t)\|_\beta^2 \lesssim  \int_0^\infty |H_1(k,t)|^2 \int_0^\infty \frac{(1 - e^{-k z})^2}{z^{1+2\beta}} \, dz \, dk \lesssim \int_0^\infty  |\widehat{{h^\dagger}}(k^2)|^2 k^{2\beta} \, dk,
\end{equation*}
with the last inequality due to the definition (A.26) of $H_1$ and equation \ref{eq:2.9}.
Then, making the change of variable $\tau = k^2$, we obtain
\begin{multline}
\|v_1(t)\|_\beta^2 \lesssim \int_0^\infty |\tau|^{\frac{2\beta-1}{2}} |\widehat{{h^\dagger}}(\tau)|^2 \, d\tau = \\ \int_0^1 |\tau|^{\frac{2\beta-1}{2}} |\widehat{{h^\dagger}}(\tau)|^2 \, d\tau +\int_1^\infty |\tau|^{\frac{2\beta-1}{2}} |\widehat{{h^\dagger}}(\tau)|^2 \, d\tau 
, \quad 0 < \beta < 1, \label{eq:2.10}
\end{multline}
Since $\mathrm{supp}({h^\dagger}) \subset (0, 2)$,
the Cauchy--Schwarz inequality implies
\[
\abs{\widehat{{h^\dagger}}(-k^2)}^2 \leq 
\left( \int_0^2 \abs{e^{-k^2 t}}^2\,\dd t \right)
\left( \int_0^2 \abs{{h^\dagger}(t)}^2\,\dd t \right)
\leq 2\norm{h}_{L^2_t(\R)}^2. 
\]
Therefore,
\begin{multline}
    \int_0^1 |\tau|^{\frac{2\beta-1}{2}} |\widehat{{h^\dagger}}(\tau)|^2 \, d\tau \leq 2\norm{{h^\dagger}}_{L^2_t(\R)}^2 \int_0^1 |\tau|^{\frac{2\beta-1}{2}} d\tau \lesssim \|{h^\dagger}\|_{H^{\frac{2s-1}{4}}_t(\mathbb{R})}, \quad 0<\beta<1,\quad s>\frac{1}{2}. \label{eq:A.31}
\end{multline}
\begin{equation}
    \int_1^\infty |\tau|^{\frac{2\beta-1}{2}} |\widehat{{h^\dagger}}(\tau)|^2 \, d\tau  \lesssim \|{h^\dagger}\|_{H^{\frac{2s-1}{4}}_t(\mathbb{R})}, \quad 0<\beta<1. \label{eq:A.32}
\end{equation}
Combining (A.31) and (A.32), we are done.

\textbf{(ii) The case $\beta = 0$.} Differentiating formula (2.11) and taking its $L^2$-norm, we have
\[
\|\partial_x^j v_1(t)\|_{L^2_x(0,\infty)} \leq \left\| \int_0^\infty e^{-\frac{\sqrt{2}}{2} k x} k^{j} |\widehat{{h^\dagger}}(k^2)| \, dk \right\|_{L^2_x(0,\infty)}.
\]
Therefore, using again Lemma 2.2 for the $L^2$-boundedness of the Laplace transform and estimating as in case (i), we find
\[
\|\partial_x^j v_1(t)\|_{L^2_x(0,\infty)} \lesssim \|k^{j} \widehat{{h^\dagger}}(k^2)\|_{L^2_k(0,\infty)} \lesssim \|{h^\dagger}\|_{H^{\frac{2j-1}{4}}_t(\mathbb{R})}, \quad j = 1,2,\dots,\lfloor s \rfloor,
\]
so that
\begin{equation}
\|v_1(t)\|_{H^{\lfloor s' \rfloor}_x(0,\infty)} \lesssim \|{h^\dagger}\|_{H^{\frac{2\lfloor s' \rfloor - 1}{4}}_t(\mathbb{R})}, \quad \lfloor s' \rfloor \geq 0. \label{eq:2.11}
\end{equation}
For $j=0$, 
\begin{equation}
    \|v_1(t)\|_{L^2_x(0,\infty)}^2 \lesssim \| \widehat{{h^\dagger}}(k^2)\|_{L^2_k(0,\infty)}^2=\int_0^1 |\widehat{{h^\dagger}}(k^2)| dk +\int_1^\infty |\widehat{{h^\dagger}}(k^2)| dk . \label{eq:A.34}
\end{equation}
As above,
\begin{equation}
    \int_0^1 |\widehat{{h^\dagger}}(k^2)| dk \lesssim \|{h^\dagger}\|_{H^{\frac{2s-1}{4}}_t(\mathbb{R})} \text{(Similar to (A.31))}. \label{eq:A.35}
\end{equation}
\begin{equation}
    \int_1^\infty |\widehat{{h^\dagger}}(k^2)| dk \lesssim \|{h^\dagger}\|_{H^{\frac{2s-1}{4}}_t(\mathbb{R})} \text{(as usual analysis)}. \label{eq:A.36}
\end{equation}

\noindent \textbf{(iii) The case $\lfloor s' \rfloor \ne 0$ and $\beta \ne 0$.} We now have $s'= \lfloor s' \rfloor + \beta$ with $\lfloor s' \rfloor \geq 1$ and $0 < \beta < 1$. Thanks to our earlier work in cases (i) and (ii), it suffices to estimate the fractional norm $\|\partial_x^{\lfloor s '\rfloor} v_1(t)\|_\beta$. Proceeding as in case (i), separating low and high frequencies, we obtain
\[
\|\partial_x^{\lfloor s \rfloor} v_1(t)\|_\beta^2 \lesssim \int_0^\infty k^{2(\lfloor s \rfloor + \beta) } |\widehat{{h^\dagger}}(k^2)|^2 \, dk = \|k^{\lfloor s \rfloor + \beta } \widehat{{h^\dagger}}(k^2)\|_{L^2_k(0,\infty)}^2,
\]
which can be estimated like the corresponding term in case (ii) to yield
\begin{equation}
\|\partial_x^{\lfloor s' \rfloor} v_1(t)\|_\beta \lesssim \|{h^\dagger}\|_{H^{\frac{2s'-1}{4}}_t(\mathbb{R})}, \quad \lfloor s \rfloor \geq 1. \label{eq:2.12}
\end{equation}
Estimates~\eqref{eq:2.10}--\eqref{eq:2.12} combined with the definition~\eqref{eq:2.7} imply the space estimate (2.17) for $v_1$. As noted earlier, this estimate can be established for $v_2$ in the exact same way. Hence, the proof of the space estimate (2.17) for $v$ is complete.

The time estimate (2.18) is easier to obtain. In particular, the change of variable $k = \sqrt{\tau}$ turns formula (2.3) into
\[
v_x(x, t) \simeq \int_0^\infty e^{i a^3 \sqrt{\tau} x + i \tau t} \widehat{{h^\dagger}}(\tau) \, d\tau + \int_{-\infty}^0 e^{i a \sqrt{-\tau} x + i \tau t} \widehat{{h^\dagger}}(\tau) \, d\tau.
\]
Therefore, since the $x$-exponentials involved above are bounded by 1 for $x \geq 0$, we find
\begin{multline}
    \|v(x)\|_{H^{\frac{2{s'}-1}{4}}_t(\mathbb{R})}^2 \lesssim \int_0^\infty (1 + \tau^2)^{\frac{2{s'}-1}{4}} e^{-\sqrt{2} \sqrt{\tau} x} |\widehat{{h^\dagger}}(\tau)|^2 \, d\tau \\ + \int_{-\infty}^0 (1 + \tau^2)^{\frac{2{s'}-1}{4}} e^{-\sqrt{2} \sqrt{-\tau} x} |\widehat{{h^\dagger}}(\tau)|^2 \, d\tau \leq \|{h^\dagger}\|_{H^{\frac{2{s'}-1}{4}}_t(\mathbb{R})}^2. \label{eq:A.38}
\end{multline}
To prove estimate (2.19), we need to use physical space decomposition:
\begin{equation}
     \|v(x)\|_{H^{\frac{2{s'}+1}{4}}_t(0,2)}^2  \simeq\|v(x)\|_{L^2_t(0,2)}^2 +  \|v_t(x)\|_{H^{\frac{2{s'}-3}{4}}_t(0,2)}^2. \label{eq:A.39}
\end{equation}
Estimating the second component on the RHS of (A.39) is routine.
\begin{multline}
    \|v_t(x)\|_{H^{\frac{2{s'}-3}{4}}_t(\mathbb{R})}^2 \lesssim \int_0^\infty (1 + \tau^2)^{\frac{2{s'}-3}{4}} \sqrt{\tau} e^{-\sqrt{2} \sqrt{\tau} x} |\widehat{{h^\dagger}}(\tau)|^2 \, d\tau \\ + \int_{-\infty}^0 (1 + \tau^2)^{\frac{2{s'}-3}{4}}\sqrt{\tau}  e^{-\sqrt{2} \sqrt{-\tau} x} |\widehat{{h^\dagger}}(\tau)|^2 \, d\tau \leq \|{h^\dagger}\|_{H^{\frac{2{s'}-1}{4}}_t(\mathbb{R})}^2. \label{eq:A.40}
\end{multline}
To estimate $ \|v(x)\|_{L^2_t(0,2)}$, we need to resort to the space estimate (2.17), which has just been proven, and the fact that $s'>\frac{1}{2}$: 
\begin{multline}
    \|v(x)\|_{L^2_t(0,2)}^2 \leq 2 \|v(x)\|_{L^\infty_t(0,2)}^2 \leq  2 \sup_{t\in (0,2)} \|v(t)\|_{L^\infty_x(\R^+)}^2 \\ \lesssim  \sup_{t\in (0,2)} \|v(t)\|_{H^{{s'}}_x(\mathbb{R}^+)} \lesssim \|{h^\dagger}\|_{H^{\frac{2{s'}-1}{4}}_t(\mathbb{R})}^2, \quad {s'}> \frac{1}{2}. \label{eq:A.41}
\end{multline}
Combining all the estimates (A.38)-(A.41), we get estimate (2.19). Thus, the proof of Theorem 2.3 is complete. Combining Theorem 2.2 and Theorem 2.3, we get Theorem 2.1.
$\square$

\begin{proof}[Proof of Lemma 2.4]
    Examining $v(x,t)$ and $\overline{v(x,t)}$, it is not hard to establish that $v(x,t)$ is real-valued.

    \noindent We claim that in the open domain $(x,t)$, $x>0$, $0<t\leq T$, the heat equation $v_t=v_{xx}$ is satisfied in the classical sense. We only need to prove this in the region $x\geq x_0$, $T\geq t\geq t_0$, for any $x_0>0$, $t_0>0$. The general result follows as $x_0$, $t_0$ approach zero. 

    We use the decomposition $v(x,t)= S[V_0](x,t)+S[0, h^\dagger](x,t)$. We want to prove high $H^s$ regularity by using exponential decay in the region $x\geq x_0$, $T\geq t\geq t_0$, then apply Sobolev embedding theorems.
    
Firstly, for $V(x,t)=S[V_0](x,t)$. For space estimate,
\begin{align}
    \norm{S[V_0](t) }_{H^3_x (\R)} =\norm{e^{-k^2t}\langle k\rangle^3 \widehat{V}_0(k) }_{L_2} \lesssim \norm{ e^{-k^2t}\langle k\rangle^3}_{L^\infty} \norm{\widehat{V}_0}_{L_2} \lesssim_{t_0} \norm{V_0}_{H^{s'}}. \label{eq:A.42}
\end{align}
For time estimate, similarly to equations (A.4-A.6), we have
\begin{multline}
    \|V_{tt}(x)\|_{L^2_t(t_0,T)}^2 \approx \int_{t_0}^T \left[ \int_0^\infty  e^{-\xi^2 t} \xi^4 (|\widehat{V}_0(-\xi)| + |\widehat{V}_0(\xi)|) \, d\xi \right]^2 dt \\
    \lesssim \norm{e^{-\xi^2 t_0} \xi^4}_{L^\infty_\xi}^2 (T-t_0) \norm{V_0}^2_{H^{s'}}. \label{eq:A.43}
\end{multline}
    Thus, combining with the restriction of equations A.11 and A.12, we have
   \begin{equation}
       \norm{ V(x)}_{H^2_t (t_0,T)} \lesssim_{t_0} \sqrt{T-t_0} \norm{V_0}_{H^{s'}}. \label{eq:A.44}
   \end{equation}

Secondly, for $S[0,h^\dagger](x,t)=v_1(x,t)+v_2(x,t)$, for $j=1,2,3$, similar to equation A.33, we have 
\begin{multline}
\|\partial_x^j v_1(t)\|_{L^2_x(x_0,\infty)} \leq \left\| \int_0^\infty e^{-\frac{\sqrt{2}}{2} k x} k^{j} |\widehat{{h^\dagger}}(k^2)| \, dk \right\|_{L^2_x(x_0,\infty)} \\ \lesssim \norm{e^{-\frac{\sqrt{2}}{4} k x_0} k^{j} }_{L^\infty_k (0,\infty)} \norm{h^\dagger}_{H^{\frac{2s-1}{4}}_t (\R)}. \label{eq:A.45}
\end{multline}
Similar estimate holds for $v_2$. Combining with the estimate in the case $j=0$ (the restriction of equations A.35-A.36), we deduce that 
\begin{equation}
    \| S[0,h^\dagger](t)\|_{H^3_x(x_0,\infty)} \lesssim_{x_0} \norm{e^{-\frac{\sqrt{2}}{4} k x_0} k^{j} }_{L^\infty_k (0,\infty)} \norm{h^\dagger}_{H^{\frac{2s-1}{4}}_t (\R)}. \label{eq:A.46}
\end{equation}
For the time estimate, 
\begin{multline}
    \norm{S[0,h^\dagger]_{tt}(x)}^2_{L^2_t(\R)}\lesssim \int_0^\infty \tau\sqrt{\tau} e^{-\sqrt{2} \sqrt{\tau} x} |\widehat{{h^\dagger}}(\tau)|^2 \, d\tau \\ + \int_{-\infty}^0 \tau \sqrt{\tau}  e^{-\sqrt{2} \sqrt{-\tau} x} |\widehat{{h^\dagger}}(\tau)|^2 \, d\tau \lesssim\norm{\tau\sqrt{\tau} e^{-\sqrt{2} \sqrt{\tau} x_0} }_{L^\infty_\tau (0,\infty)} \norm{ h^\dagger}_{L^2_t (\R)}\lesssim_{x_0}  \|{h^\dagger}\|_{H^{\frac{2{s'}-1}{4}}_t(\mathbb{R})}^2. \label{eq:A.47}
\end{multline}
Combining with the restriction of equations A.40 and A.41, for $x\geq x_0$, we have $S[0,h^\dagger](x) \in H^2_t (\R)$.

Thus, in the region $x_0\leq x,$ $t_0\leq t\leq T$, $v(x,t)=S[V_0]+S[0,v_0]$ is at least $H^3$ in the $x$-direction and $H^2$ in the $t$-direction. Hence, by Sobolev embedding theorems, $v(x,t)$ is twice differentiable in the $x$-direction, and once differentiable in the $t$ direction, in the classical sense.
\end{proof}
\section{Appendix B: Proof of Lemma 8.2 (Quotient Lemma)}\label{app:B}\label{sec:B}

Throughout this appendix, all functions are real
valued. Also, we need the following result.

\begin{prop}[Morrey's inequality in one dimension]\label{fact:B-morrey}
If $v\in H^{1}(\mathbb{R}^{+})$ then $v$ has an absolutely continuous
representative and, for all $x,y\in\mathbb{R}^{+}$,
\begin{equation}\label{eq:B-morrey}
|v(x)-v(y)|
=\Bigl|\int_{y}^{x}v'(t)\,dt\Bigr|
\le \Bigl(\int_{\min(x,y)}^{\max(x,y)}|v'|^{2}\Bigr)^{1/2}|x-y|^{1/2}
\le \|v'\|_{L^{2}(\mathbb{R}^{+})}\,|x-y|^{1/2},
\end{equation}
i.e. $H^{1}(\mathbb{R}^{+})\hookrightarrow C^{0,1/2}(\overline{\mathbb{R}^{+}})$.
This is Morrey's inequality $W^{1,p}\hookrightarrow C^{0,1-n/p}$ with $n=1$,
$p=2$.
\end{prop}

Also, we note the following elementary algebraic identity: if $G=1+g\ge c>0$, then for all $x,y\in\mathbb{R}^{+}$
\begin{equation}\label{eq:B-identity}
\frac{f(x)}{G(x)}-\frac{f(y)}{G(y)}
=\frac{f(x)-f(y)}{G(x)}-f(y)\,\frac{g(x)-g(y)}{G(x)G(y)} .
\end{equation}

Since $s>1/2$, and from the assumption, we have the following estimates:
\begin{align}
    \|f\|_{L^{\infty}}&\lesssim_s \|f\|_{H^{s}} \label{eq:B.3}\\
    \|g\|_{L^{\infty}}&\lesssim_s\|g\|_{H^{s}} \label{eq:B.4}\\
    \|G\|_{L^{\infty}}&\le 1+C(s)\|g\|_{H^{s}} \label{eq:B.5}\\
    \|1/G\|_{L^{\infty}}&\le c^{-1}. \label{eq:B.6}
\end{align}

The case $1<s<\tfrac32$ requires the following product estimate, which is
adapted to the fact that for $\sigma:=s-1<\tfrac12$ the functions $f',g'\in
H^{\sigma}$ might not be bounded.

\begin{lemma}\label{lem:B-product}
Let $0<\sigma<\tfrac12$. Let $v\in L^{2}(\mathbb{R}^{+})$ with
$\norm{v}_{\sigma}<\infty$, and let $w\in L^{\infty}(\mathbb{R}^{+})$ satisfy the
H\"older bound
\[
|w(x)-w(y)|\le L\,|x-y|^{1/2}\qquad\text{for all }x,y\in\mathbb{R}^{+}.
\]
Then $vw\in H^{\sigma}(\mathbb{R}^{+})$ and
\begin{equation}\label{eq:B-product}
\norm{vw}_{\sigma}
\le \|w\|_{L^{\infty}}\norm{v}_{\sigma}
+\Bigl(\frac{2}{1-2\sigma}\Bigr)^{1/2}L\,\|v\|_{L^{2}}
+\Bigl(\frac{4}{\sigma}\Bigr)^{1/2}\|w\|_{L^{\infty}}\|v\|_{L^{2}} .
\end{equation}
\end{lemma}

\begin{proof}
Let $d\mu(x,y):=|x-y|^{-1-2\sigma}dx\,dy$ on
$\mathbb{R}^{+}\times\mathbb{R}^{+}$, so that
$\norm{\cdot}_{\sigma}=\|\,\Delta(\cdot)\,\|_{L^{2}(\mu)}$ with
$\Delta v(x,y):=v(x)-v(y)$. Splitting the increment of the product,
\[
\Delta(vw)(x,y)=w(x)\,\Delta v(x,y)+v(y)\,\Delta w(x,y),
\]
the triangle inequality in $L^{2}(\mu)$ gives
\begin{equation}\label{eq:B-split}
\norm{vw}_{\sigma}\le \|w\|_{L^{\infty}}\norm{v}_{\sigma}
+\Bigl(\iint_{x,y\in\R^+} \frac{|v(y)|^{2}|w(x)-w(y)|^{2}}{|x-y|^{1+2\sigma}}\,dx\,dy\Bigr)^{1/2}.
\end{equation}
We estimate the last integral by splitting the region of integration according
to $|x-y|\le 1$ and $|x-y|>1$, and by writing $x=y+h$ for the inner variable.

 If $|x-y|\le 1$, the H\"older hypothesis gives
$|w(x)-w(y)|^{2}\le L^{2}|x-y|$, hence
\[
\iint_{|x-y|\le 1}\frac{|v(y)|^{2}|w(x)-w(y)|^{2}}{|x-y|^{1+2\sigma}}\,dx\,dy
\le L^{2}\int_{0}^{\infty}|v(y)|^{2}
\Bigl(\int_{|h|\le 1}|h|^{-2\sigma}\,dh\Bigr)dy
=\frac{2L^{2}}{1-2\sigma}\|v\|_{L^{2}}^{2},
\]
where we used $\int_{|h|\le1}|h|^{-2\sigma}dh=\frac{2}{1-2\sigma}$, which is
finite precisely because $\sigma<\tfrac12$.

 If $|x-y|>1$, as
$|w(x)-w(y)|\le 2\|w\|_{L^{\infty}}$ and
$\int_{|h|>1}|h|^{-1-2\sigma}dh=\tfrac1\sigma$,
\[
\iint_{|x-y|>1}\frac{|v(y)|^{2}|w(x)-w(y)|^{2}}{|x-y|^{1+2\sigma}}\,dx\,dy
\le \frac{4\|w\|_{L^{\infty}}^{2}}{\sigma}\,\|v\|_{L^{2}}^{2}.
\]
Inserting the two bounds into \eqref{eq:B-split} and using
$(a+b)^{1/2}\le a^{1/2}+b^{1/2}$ yields \eqref{eq:B-product}. Finally
$\|vw\|_{L^{2}}\le\|w\|_{L^{\infty}}\|v\|_{L^{2}}<\infty$, so
$vw\in H^{\sigma}(\mathbb{R}^{+})$.
\end{proof}

\begin{proof}[Proof of Lemma 8.2]
We treat separately the three ranges $\tfrac12<s<1$, $s=1$ and
$1<s<\tfrac32$. In all three cases we abbreviate
\[
F:=\|f\|_{H^{s}(\mathbb{R}^{+})},\qquad
M:=1+\|g\|_{H^{s}(\mathbb{R}^{+})} ,
\]
so that $\|f\|_{L^{\infty}}\le C(s)F$ and
$\|G\|_{L^{\infty}}\le C(s)M$, while $c\le G$ pointwise. Note also that
$u$ is measurable and $|u|\le c^{-1}|f|$, whence
\begin{equation}\label{eq:B-L2}
\|u\|_{L^{2}(\mathbb{R}^{+})}\le c^{-1}\|f\|_{L^{2}(\mathbb{R}^{+})}\le c^{-1}F .
\end{equation}

\medskip
\noindent\textbf{Case 1: $\tfrac12<s<1$.}
Since $G(x)-G(y)=g(x)-g(y)$
and $G\ge c$, the identity \eqref{eq:B-identity} yields the pointwise bound
\begin{equation}\label{eq:B-case1-pt}
|u(x)-u(y)|
\le c^{-1}|f(x)-f(y)|+c^{-2}\|f\|_{L^{\infty}}\,|g(x)-g(y)|,
\qquad x,y\in\mathbb{R}^{+} .
\end{equation}
Squaring \eqref{eq:B-case1-pt}, dividing by $|x-y|^{1+2s}$, integrating over
$\mathbb{R}^{+}\times\mathbb{R}^{+}$ and using the triangle inequality in
$L^{2}(\mu)$ (notation as in the proof of Lemma~\ref{lem:B-product}), we obtain
\[
\norm{u}_{s}\le c^{-1}\norm{f}_{s}+c^{-2}\|f\|_{L^{\infty}}\norm{g}_{s}
\le c^{-1}F+C(s)c^{-2}F\,\|g\|_{H^{s}}
\le C(s,c)\,F\,M .
\]
Together with \eqref{eq:B-L2} this gives $u\in H^{s}(\mathbb{R}^{+})$ and
$\|u\|_{H^{s}}\le C(s,c)FM$.

\medskip
\noindent\textbf{Case 2: $s=1$.}
For $u=f\cdot\frac1G$, applying Leibniz rule, we have:
\begin{equation}\label{eq:B-leibniz}
u'=\frac{f'}{G}-\frac{f\,g'}{G^{2}}\qquad\text{a.e. on }\mathbb{R}^{+}.
\end{equation}
Therefore, using $G\ge c$,
\[
\|u'\|_{L^{2}}
\le c^{-1}\|f'\|_{L^{2}}+c^{-2}\|f\|_{L^{\infty}}\|g'\|_{L^{2}}
\le c^{-1}F+C(s)c^{-2}F\|g\|_{H^{1}}
\le C(s,c)\,F\,M ,
\]
which with \eqref{eq:B-L2} gives $u\in H^{1}(\mathbb{R}^{+})$ and
$\|u\|_{H^{1}}\le C(s,c)FM$.

\medskip
\noindent\textbf{Case 3: $1<s<\tfrac32$.}
Set $\sigma:=s-1\in(0,\tfrac12)$. Since
$H^{s}(\mathbb{R}^{+})\hookrightarrow H^{1}(\mathbb{R}^{+})$, Case~2 applies and
yields $u\in H^{1}(\mathbb{R}^{+})$ with
\begin{equation}\label{eq:B-H1bound}
\|u\|_{H^{1}}\le C(s,c)\,F\,M ,
\end{equation}
together with equation \eqref{eq:B-leibniz}. Hence, it remains to bound $\norm{u'}_{\sigma}$. Write
\[
u'=A-B,\qquad
A:=f'\cdot\frac1G,
\qquad
B:=g'\cdot\frac{f}{G^{2}} ,
\]
so that $\norm{u'}_{\sigma}\le \norm{A}_{\sigma}+\norm{B}_{\sigma}$. The hypotheses
of Lemma~\ref{lem:B-product} is satisfied for each product; note that by definition, 
$f',g'\in L^{2}$ with $\norm{f'}_{\sigma}\le F$ and $\norm{g'}_{\sigma}\le\|g\|_{H^{s}}\le M$.

\smallskip
\emph{(i) The term $A=f'\cdot\frac1G$.} We have
$\|\frac1G\|_{L^{\infty}}\le c^{-1}$, and by Morrey's inequality \eqref{eq:B-morrey} applied to $g$,
\[
\left|\frac1{G(x)}-\frac1{G(y)}\right|
=\frac{|g(y)-g(x)|}{G(x)G(y)}
\le c^{-2}\|g'\|_{L^{2}}\,|x-y|^{1/2},
\]
so $\frac1G$ satisfies the hypothesis of Lemma~\ref{lem:B-product} with
$L_{1}:=c^{-2}\|g'\|_{L^{2}}\le c^{-2}M$. Hence, by \eqref{eq:B-product},
\begin{equation}\label{eq:B-A}
\norm{A}_{\sigma}
\le c^{-1}\norm{f'}_{\sigma}
+C(\sigma)\,c^{-2}M\,\|f'\|_{L^{2}}
+C(\sigma)\,c^{-1}\|f'\|_{L^{2}}
\le C(s,c)\,F\,M .
\end{equation}

\smallskip
\emph{(ii) The term $B=g'\cdot\frac{f}{G^{2}}$.} By Sobolev embedding,
$\|\frac{f}{G^{2}}\|_{L^{\infty}}\le c^{-2}\|f\|_{L^{\infty}}\le C(s)c^{-2}F$. Moreover,
for $x,y\in\mathbb{R}^{+}$,
\[
\left|\frac{f}{G^{2}}(x)-\frac{f}{G^{2}}(y)\right|
\le \frac{|f(x)-f(y)|}{G(x)^{2}}
+\|f\|_{L^{\infty}}\frac{|G(y)^{2}-G(x)^{2}|}{G(x)^{2}G(y)^{2}} ,
\]
and since
$|G(x)^{2}-G(y)^{2}|=\bigl(G(x)+G(y)\bigr)|g(x)-g(y)|
\le 2\|G\|_{L^{\infty}}\|g'\|_{L^{2}}|x-y|^{1/2}$
by \eqref{eq:B-morrey}, we obtain
\[
\left|\frac{f}{G^{2}}(x)-\frac{f}{G^{2}}(y)\right|\le L_{2}\,|x-y|^{1/2},
\qquad
L_{2}:=c^{-2}\|f'\|_{L^{2}}
+2c^{-4}\|G\|_{L^{\infty}}\|f\|_{L^{\infty}}\|g'\|_{L^{2}} .
\]
Notice that $L_{2}\le C(s,c)\,F\,M^{2}$. Lemma~\ref{lem:B-product}
therefore gives
\begin{equation}\label{eq:B-B}
\norm{A}_{\sigma}
\le \left \|\frac{f}{G^{2}}\right\|_{L^{\infty}}\norm{g'}_{\sigma}
+C(\sigma)\bigg(L_{2}+\left \|\frac{f}{G^{2}}\right\|_{L^{\infty}}\bigg)\|g'\|_{L^{2}}
\le C(s,c)\,F\,M^{3}.
\end{equation}

\smallskip
Combining \eqref{eq:B-A}, \eqref{eq:B-B}, and \eqref{eq:B-H1bound}, we conclude that $u\in H^{s}(\mathbb{R}^{+})$ with
\[
\|u\|_{H^{s}(\mathbb{R}^{+})}
\le \|u\|_{H^{1}}+\norm{A}_{\sigma}+\norm{B}_{\sigma}
\le C(s,c)\,\|f\|_{H^{s}(\mathbb{R}^{+})}
\bigl(1+\|g\|_{H^{s}(\mathbb{R}^{+})}\bigr)^{3}.
\]
This completes the proof in all three cases.
\end{proof}

\begin{remark}\label{rem:B-sharp}
The exponent $3$ in Lemma 8.2 is not optimal; it arises from
the crude bound $L_{2}\lesssim FM^{2}$ in step (ii) of Case~3. In the range
$\tfrac12<s\le1$ the proof gives the linear dependence
$\|u\|_{H^{s}}\le C(s,c)\|f\|_{H^{s}}(1+\|g\|_{H^{s}})$.
\end{remark}

\bibliographystyle{elsarticle-num} 
\bibliography{main}
\end{document}